\documentclass[twocolumn]{autart}    
\newcommand{\tr}{\mathbf{tr}}
\newcommand{\E}{\mathbf{E}}

\usepackage{graphicx,import}          
\usepackage[french,english]{babel}
\usepackage{color}
\usepackage{amsmath,amssymb,amsfonts}  
\usepackage{booktabs}
\usepackage[table,xcdraw]{xcolor}
\everymath{\displaystyle}
\usepackage{hyperref}
\usepackage[authoryear]{natbib}
\newtheorem{assumption}{Assumption}

\newtheorem{lemma}{Lemma}

\begin{document}
\begin{frontmatter}
\title{A Nonlinear Navigation Observer Using IMU and Generic Position Information}
\thanks[footnoteinfo]{This work was supported by the National Sciences and Engineering Research Council of Canada (NSERC), under the grants NSERC-DG 228465-2013-RGPIN and NSERC-DG RGPIN-2020-04759. A preliminary and partial version of this work was presented in \citep{berkane2019position}.}
\author[Soulaimane]{Soulaimane Berkane}\ead{soulaimane.berkane@uqo.ca},               
\author[hamid]{Abdelhamid Tayebi}\ead{atayebi@lakeheadu.ca},
\author[simone]{Simone de Marco}\ead{sdemarco@i3s.unice.fr}

\address[Soulaimane]{ Department of Computer Science and Engineering, University of Quebec in   Outaouais,   Gatineau, Quebec,   Canada}
\address[hamid]{Department of Electrical Engineering, Lakehead University, Thunder Bay, Ontario, Canada.}
\address[simone]{I3S, Universit\'e de Nice, Sophia-Antipolis, CNRS, France.}


\selectlanguage{english}
\begin{abstract}
This paper deals with the problem of full state estimation for vehicles navigating in a three dimensional space. We assume that the vehicle is equipped with an Inertial Measurement Unit (IMU) providing body-frame measurements of the angular velocity, the specific force, and the Earth's magnetic field. Moreover, we consider available sensors that provide partial or full information about the position of the vehicle. Examples of such sensors are those which provide full position measurements (\textit{e.g.,} GPS), range measurements (\textit{e.g.,} Ultra-Wide Band (UWB) sensors),  inertial-frame bearing measurements (\textit{e.g.,} motion capture cameras), or altitude measurements (altimeter). We propose a generic semi-globally exponentially stable nonlinear observer that estimates the position, linear velocity, linear acceleration, and attitude of the vehicle, as well as the gyro bias. We also provide a detailed observability analysis for different types of measurements. Simulation and experimental results are provided to demonstrate the effectiveness of the proposed estimation scheme.
\end{abstract}
\end{frontmatter}
\section{Introduction}
\subsection{Motivation}
Inertial navigation systems (INS) are of great importance in many autonomous vehicles and robot platforms \citep{Titterton2004StrapdownTechnology}. They combine measurements from translational motion sensors (accelerometers) and rotational motion sensors (gyroscopes), to track the position, velocity and orientation of a vehicle with respect to a reference frame. Three orthogonal rate gyroscopes and three orthogonal accelerometers, measuring angular velocity and linear acceleration respectively, are typically included in an Inertial Measurement Unit (IMU) which is used, in addition to a processing unit, inside an INS. However, the use of an INS alone for navigation usually leads to unreliable state estimates since measurement errors and unknown initial conditions cause the estimation error to drift over time \citep{Woodman2007AnNavigation}. In this paper, we tackle the problem of the design of inertial navigation systems aided with different types of position information. The approach we take consists in using a nonlinear observer that fuses IMU and position sensors to provide reliable state estimation with semi-global exponential stability guarantees.
\subsection{Prior Works}
Inertial navigation systems are usually aided by position sensors such as Global Positioning System (GPS) which allows to correct the position estimates over time, thus, keeping the estimation errors small and bounded; see, {\it e.g.,} \citep{vik2001nonlinear,Farrell2008AidedSensors,Grip2013}. Other type of sensors that can provide range (distance)  measurements  to  known  source  points can also be used to provide position information. For instance, there are some GPS receivers that provide access to the raw GPS observations (pseudo-range measurements) which are used in {\it tightly} coupled GPS/INS integration schemes. An advantage of the tightly coupled integration over the \textit{loosely} coupled integration is the possibility to use few raw GPS pseudo-ranges that would otherwise be insufficient to provide full position estimates \citep{johansen2015nonlinear,Johansen2017NonlinearMeasurements,Bryne2017NonlinearAspects}.  Another type of range measurement sensors that are getting popular in indoor applications is the Ultra Wide-Band (UWB) radio technology \citep{Gryte2017RobustGPS,Hamer2018Self-calibratingLocalization}. The idea consists in mounting a network of radio modules (anchors) at known locations along with a receiver on the vehicle. By communicating signals between the anchors and the receiver, the vehicle is able to calculate the distance (range) to the transmitting anchor in the same way a GPS receiver communicates with the satellites. The UWB-based localization technology has shown very promising accuracy in short-range applications. Utra-short/Long baseline (USBL/LBL) sensors are other examples of range sensors used in marine applications  \citep{Batista2011SingleDesign,Batista2016TightlySystem}. Finally, raw data from motion capture systems ({\it e.g.,} OptiTrack, Vicon, Xsens) can also be used as a source of position information providing inertial-frame bearing measurements representing the projections of the relative position vectors (with respect to the fixed cameras) on the unit sphere \citep{Batista2015NavigationMeasurements,Batista2013GloballyMeasurements,Hamel2017PositionMeasurements}.  

Although Kalman-type filters, such as \citep{sabatini2006quaternion,Farrell2008AidedSensors,Crassidis2006Sigma-pointNavigation,Whittaker2017InertialRepresentations}, are considered industry-standard solutions for inertial navigation systems, these stochastic filters are often based on linearization assumptions and may fail when the initial estimation errors are large.  Recently, nonlinear state observers, also called deterministic estimators, have been developed for autonomous navigation applications.  Attitude observers have been proposed in \citep{ Mahony2008NonlinearGroup,Batista2012GloballyEstimation,Zlotnik2016NonlinearDirectly,Berkane2016GlobalSystems,Berkane2017OnSO3} to name few. Velocity-aided attitude observers can be found for instance in \citep{Bonnabel2006,hua2010attitude,roberts2011,Berkane2017AttitudeMeasurements}. Full state observers (attitude, position and linear velocity) are developed in \citep{Johansen2017NonlinearMeasurements,hansen2018nonlinear, Bryne2017NonlinearAspects}. The advantage of the nonlinear observers is their theoretically proven stability guarantees, as well as their computational simplicity compared to the stochastic filters. 
\subsection{Contributions and Paper Organization}

In this work, we propose a nonlinear observer, relying on inertial measurements and full or partial position information, for the simultaneous estimation of the position, velocity, acceleration, attitude, and gyro bias of a rigid body system. The proposed observer achieves semi-global exponential stability, which is, to the best of the authors knowledge, the strongest stability result available in the literature dealing with the problem at hand.
The attitude estimates are directly obtained on the Special Orthogonal group of rotations $\mathbb{SO}(3)$, thus avoiding any singularities or ambiguities related to the use of other attitude parametrizations. The attitude estimation part is based on the nonlinear complimentary filter \citep{Mahony2008NonlinearGroup} where the innovation term is generated using the IMU measurements and the estimates of the unknown inertial linear acceleration as done, for instance, in \citep{hua2010attitude,roberts2011,Grip2013,Berkane2017AttitudeMeasurements,Johansen2017NonlinearMeasurements}. The translational motion observer is designed to handle, in a unified manner, different types of position sensors; a feature that cannot be found in the existing observers in the literature that are usually tailored to the type of sensors used. The translational estimation part is a high-gain observer with a similar structure to the nonlinear observers in \citep{Grip2013,Bryne2017NonlinearAspects,Johansen2017NonlinearMeasurements}. Another important contribution of this work consists in the detailed uniform observability analysis carried out for the different types of sensors used.
Finally, we would like to point out that a preliminary version of this work appeared in our conference paper \citep{berkane2019position}. In this extended version, we 1) provide complete proofs of the main results, 2) derive a generic observability condition on the output matrix, and conduct a detailed observability analysis with different types of position measurements, and 3) experimentally validate the proposed observer on a quadrotor UAV.

The paper is organized as follows. After some preliminaries in Section \ref{section:preliminaries}, we formulate our estimation problem in Section \ref{section:problem} where we give details about the considered vehicle's model, the possible available measurements and the technical assumptions needed for our main result. Then, in Section \ref{section:riccati_observer}, the proposed nonlinear observer is provided and the main result is announced in Section \ref{section:stability}. In Section \ref{section:observability} we study the observability of the translational motion for different scenarios of position measurements. Simulation and experimental results are provided in Sections \ref{section:simulation} and \ref{section:experiments}, respectively. Finally, Section \ref{section:conclusion} wraps up the paper with concluding remarks. The appendices are devoted for the technical proofs of our theoretical results.

%
\section{Preliminaries}\label{section:preliminaries}
We denote by $\mathbb{R}$ the set of reals and by $\mathbb{N}$ the set of natural numbers. We denote by $\mathbb{R}^n$ the $n$-dimensional Euclidean space, by $\mathbb{S}^n$ the unit $n$-sphere embedded in $\mathbb{R}^{n+1}$ and by $\mathbb{B}_\epsilon=\{x\in\mathbb{R}^3: \|x\|\leq \epsilon\}$ the closed ball in $\mathbb{R}^3$ with radius $\epsilon$. We use $\|x\|$ to denote the Euclidean norm of a vector $x\in\mathbb{R}^n$ and $\|A\|_F$ to denote the Frobenius norm of a matrix $A\in\mathbb{R}^{n\times n}$. Let $I_n$ be the $n$-by-$n$ identity matrix and let $e_i$ denote the $i-$th column of $I_n$. The Special Orthogonal group of order three is denoted by $\mathbb{SO}(3) := \{ A \in \mathbb{R}^{3\times 3}:\;\mathbf{det}(A)=1,\; AA^{\top}=A^\top A= I_3 \}$. The set $\mathfrak{so}(3):=\left\{\Omega\in\mathbb{R}^{3\times 3}\mid\;\Omega^{\top}=-\Omega\right\}$ denotes the Lie algebra of $\mathbb{SO}(3)$. For $x,~y\in\mathbb{R}^3$, the map $[\cdot]_\times: \mathbb{R}^3\to\mathfrak{so}(3)$ is defined such that $[x]_\times y=x\times y$ where $\times$ is the vector cross-product on $\mathbb{R}^3$. The inverse isomorphism of the map $[\cdot]_\times$ is defined by $\mathrm{vex}:\mathfrak{so}(3)\to\mathbb{R}^3$, such that $\mathrm{vex}([\omega]_\times)=\omega$, for all $\omega\in\mathbb{R}^3$ and $[\mathrm{vex}(\Omega)]_\times=\Omega,$ for all $\Omega\in\mathfrak{so}(3)$. The composition map $\psi := \mathrm{vex}\circ \mathbf{P}_{\mathfrak{so}(3)}$ extends the definition of
        $\mathrm{vex}$ to  $\mathbb{R}^{3\times 3}$, where $\mathbf{P}_{\mathfrak{so}(3)}:\mathbb{R}^{3\times 3}\to\mathfrak{so}(3)$ is the projection map on the Lie algebra $\mathfrak{so}(3)$ such that $\mathbf{P}_{\mathfrak{so}(3)}(A):=(A-A^{\top})/2$. Accordingly, for a given $3$-by-$3$ matrix $A=[a_{ij}]_{i,j = 1,2,3}$, one has ${\textstyle\psi(A)=\frac{1}{2}[a_{32}-a_{23},a_{13}-a_{31},a_{21}-a_{12}]
}$. We define ${\textstyle|R|:=\frac{1}{4}\tr(I_3-R)=\frac{1}{8}\|I_3-R\|_F^2\in[0, 1]}$ as the normalized Euclidean distance on $\mathbb{SO}(3)$. Given a scalar $c>0$, we define the saturation function $\mathbf{sat}_c:\mathbb{R}^n\to\mathbb{R}^n$ such that
$
\mathbf{sat}_c(x):=\min(1,c/\|x\|)x.
$
Given two scalars $c,\epsilon>0$, we also define the smooth projection function $\mathbf{P}_c^\epsilon:\mathbb{R}^3\times\mathbb{R}^3\to\mathbb{R}^3$, found for instance in \citep{Krstic1995NonlinearDesign}, as follows:
\begin{equation}\label{eq:pprojection}
\mathbf{P}_c^\epsilon(\hat\phi,\mu):=
\begin{cases}
\mu,&\textrm{if}\;\|\hat\phi\|<c\;\textrm{or}\;\hat\phi^\top\mu\leq 0,\\
{\scriptstyle
\left(I-\theta(\hat\phi)\frac{\hat\phi\hat\phi^\top}{\|\hat\phi\|^2}\right)}\mu,&\textrm{otherwise},
\end{cases}
\end{equation}
where we let $\theta(\hat\phi):=\min(1,(\|\hat\phi\|-c)/\epsilon)$.  The projection operator $\mathbf{P}_c^\epsilon(\hat\phi,\mu)$ is locally Lipschitz in its arguments. Moreover, provided that $\|\phi\|\leq c$, the projection map $\mathbf{P}_c^\epsilon(\hat\phi,\mu)$ satisfies, along the trajectories of $\dot{\hat{\phi}}=\mathbf{P}_c^\epsilon(\hat\phi,\mu), \|\hat\phi(0)\|\leq c+\epsilon$, the following properties:
            \begin{align}
            \label{ineq:P1}
            &\|\hat\phi(t)\|\leq c+\epsilon,\;\forall t\geq 0,\\
            \label{ineq:P2}
            &(\hat \phi-\phi)^\top\mathbf{P}_c^\epsilon(\hat\phi,\mu)\leq(\hat\phi-\phi)^\top\mu,\\
            \label{ineq:P3}
            &\|\mathbf{P}_c^\epsilon(\hat\phi,\mu)\|\leq\|\mu\|.
           \end{align}
Finally, the pair $(A(\cdot),C(\cdot))$ is uniformly observable if there exist $\delta,\mu>0$ such that, for all $t\geq 0$, the observability Gramian matrix satisfies:
\begin{equation}\label{obs:gramian}
W(t,t+\delta):=\int_t^{t+\delta}\Phi^\top(s,t)C^\top(s)C(s)\Phi(s,t)ds\ge \mu I_n, 
\end{equation}
where $\Phi(t,s)$ is the state transition matrix associated to $A(t)\in\mathbb{R}^{n\times n}$, which is defined by $\dot\Phi(t,s)=A(t)\Phi(t,s)$ and $\Phi(t,t)=I_n$. Different (explicit) sufficient conditions for uniform observability have been developed in the literature; see, for instance, \citep{Bristeau2010DesignObservability,Hamel2017PositionMeasurements,Batista2017RelaxedRealizations,morin2017uniform}.

\section{Problem Formulation}\label{section:problem}
In this paper, we consider the following 3D kinematics of a rigid body:
\begin{align}
\label{eq:dp}
\dot p&=v,\\
\label{eq:dv}
\dot v&=ge_3+Ra_B,\\
\label{eq:dR}
\dot R&=R[\omega]_\times,
\end{align}
where $p\in\mathbb{R}^3$ is the inertial position of the vehicle's center of gravity, $v\in\mathbb{R}^3$ represents the inertial linear velocity, $R\in \mathbb{SO}(3)$ is the attitude matrix describing the orientation of a body-attached frame with respect to the inertial frame, $\omega$ is the angular velocity of the body-attached frame with respect to the inertial frame expressed in the body-attached frame, $g$ is the norm of the acceleration due to gravity, and $a_B$ is the {\it specific force}, capturing all non-gravitational forces applied to the vehicle, expressed in the body-attached frame. We may also use the term {\it apparent acceleration} interchangeably with specific force.

We assume available an  IMU that provides measurements of the angular velocity, the body-attached frame non-gravitational acceleration and the body-attached frame magnetic field. These sensors are modelled as:
\begin{align}
\label{eq:wy}
\omega^y&=\omega+b_\omega,\\
\label{eq:aB}
a_B&=R^\top a_I,\\
\label{eq:bm}
m_B&=R^\top m_I,
\end{align}
where $b_\omega$ is a constant unknown gyro bias, $m_I$ is the constant and known earth's magnetic field and $a_I(t)$ is the unknown time-varying apparent acceleration. Note that, in practice, the above sensor measurements are usually corrupted by random noise as well. However, since in this paper we are interested to develop a nonlinear deterministic observer, we have assumed that there is no stochastic noise in the measurements although deterministic observers are shown to have a good level of noise filtering in practice. For the rotational motion subsystem, the following is a general observability assumption used in the field of attitude estimation.
\begin{assumption}[Rotational Motion Observability]\label{assumption::obsv}
There exists $c_0>0$ such that $\|m_I\times a_I(t)\|\geq c_0, \forall t\geq 0$.
\end{assumption}
Assumption \ref{assumption::obsv} is guaranteed if the time-varying apparent acceleration $a_I(t)$ is non-vanishing and is \textit{always} non-collinear with the constant magnetic field vector $m_I$. Note that $a_I(t)=0$ corresponds to  the rigid body being in a free-fall case ($\dot v=ge_3$) which is not likely under normal flight conditions. 

We also assume in this work that we have measurements of the following position output vector:
\begin{align}\label{eq:yp}
y=C_p(t)p,
\end{align}
where $C_p(t)\in\mathbb{R}^{m\times 3}, m\in\mathbb{N},$ is a time-varying bounded output matrix. The measurement $y$ can be obtained from different possible sensors, depending on the application at hand, that provide some information about the position. In Section \ref{section:observability} we will discuss some particular examples of measurements that can be written as in \eqref{eq:yp}. 
Considering the extended state $x:=[p^\top,v^\top,a_I^\top]^\top\in\mathbb{R}^9$ for the translational motion, in view of \eqref{eq:dp}-\eqref{eq:dv} and \eqref{eq:yp}, the dynamics of $x$ are given by
\begin{align}
\label{eq:dx}
\dot x&=Ax+B_1ge_3+B_2\dot a_I,\\
\label{eq:y}
y&=C(t)x,
\end{align}
where $A, B_1, B_2$ and $C(t)$ are defined as follows:
\begin{align}
&A=\begin{bmatrix}
0_{3\times 3}&I_3&0_{3\times 3}\\
0_{3\times 3}&0_{3\times 3}&I_3\\
0_{3\times 3}&0_{3\times 3}&0_{3\times 3}
\end{bmatrix},
\label{eq:C}
C(t)=\begin{bmatrix}
C_p(t)^\top\\
0_{3\times m}\\
0_{3\times m}
\end{bmatrix}^\top,\\
&B_1:=\begin{bmatrix}
0_{3\times 3}\\
I_3\\
0_{3\times 3}
\end{bmatrix},
B_2:=\begin{bmatrix}
0_{3\times 3}\\
0_{3\times 3}\\
I_3
\end{bmatrix}.
\end{align}
The translational system \eqref{eq:dx}-\eqref{eq:y} is a linear time-varying system with an unknown input $\dot a_I$ (the jerk). The latter corresponds to the derivative of the apparent acceleration $a_I$ which is known only in body-frame, {\it i.e.,} $a_B=R^\top a_I$. Therefore, there is a coupling between the translational dynamics \eqref{eq:dx}-\eqref{eq:y} and the rotational dynamics \eqref{eq:dR} through the measurement equation \eqref{eq:aB} of the accelerometer. Most adhoc methods in practice assume that $a_I\approx -ge_3$ to remove this coupling between the translational and rotational dynamics. However, this assumption holds only for non-accelerated vehicles, {\it i.e., } when $\dot v\approx 0$. In this work, we instead design our estimation algorithm without this latter assumption. Therefore, our proposed approach will be most suitable for accelerated vehicles where the performance of adhoc approaches is compromised.

The objective is to design a full navigation observer that takes the measurements \eqref{eq:wy}-\eqref{eq:yp} and outputs reliable estimates for the position $p$, velocity $v$, orientation $R$, apparent acceleration $a_I$, and gyro bias $b_\omega$. More specifically, we want to design an exponentially convergent nonlinear observer that estimates the state variables $(p,v,a_I,R,b_\omega)$ under the above observability conditions and the following mild constraints on the trajectory of the vehicle:
\begin{assumption}\label{assumption::bounded_ra}
There exist constants $c_1,c_2,c_3>0$ such that $c_1\leq\|a_I(t)\|\leq c_2$ and $\|\dot a_I(t)\|\leq c_3$ for all $t\geq 0$.
\end{assumption}
\begin{assumption}\label{assumption::bounded_bw}
There exists constants $c_4,c_5>0$ such that $\|\omega(t)\|\leq c_4$ and $\|b_\omega\|\leq c_5$ for all $t\geq 0$.
\end{assumption}
Assumptions \ref{assumption::bounded_ra} and \ref{assumption::bounded_bw} impose some realistic constraints on the systems trajectory which are needed to carry out the stability analysis.

\section{Main Results}
\subsection{Observer Design}\label{section:riccati_observer}
In this section, we design our navigation observer that estimates the vehicle's state $(p,v,R)$ as well as the constant gyro bias $b_\omega$ and the unknown apparent acceleration $a_I$. We propose the following nonlinear navigation observer:
\begin{align}\label{eq:dxhat}
\hat x&=\hat z+B_2\hat Ra_B,\\
\label{eq:dzhat}
\dot{\hat z}&=A\hat x+B_1ge_3+K(t)(y-C(t)\hat x)+\sigma_1,\\
\label{eq:dRhat}
\dot{\hat R}&=\hat R[\omega^y-\hat b_\omega+k_1\sigma_2]_\times,\\
\label{eq:dbhat}
\dot{\hat b}_\omega&=\mathbf{P}_{c_5}^{\epsilon_b}(\hat b_\omega,-k_2\sigma_2),
\end{align}
with initial conditions $\hat x(0)\in\mathbb{R}^6, \hat R(0)\in\mathbb{SO}(3)$ and $\hat b_\omega(0)\in\mathbb{B}_{c_5+\epsilon_b}$. The innovation terms  $\sigma_1\in\mathbb{R}^9$  and $\sigma_2\in\mathbb{R}^3$ are defined as follows:
\begin{align}\label{eq:sigma2:1}
&\sigma_1=-k_1B_2\hat R [\sigma_2]_\times a_B,\\
&\sigma_2=\rho_1(m_B\times\hat R^\top m_I)+\rho_2(a_B\times\hat R^\top\mathbf{sat}_{\hat c_2}(B_2^\top \hat x)).
\label{eq:sigma2:2}
\end{align}
The scalars $k_1,k_2,\rho_1,\rho_2,\epsilon_b, \hat c_2$ are positive tuning parameters with $\hat c_2>c_2$, the parameters $c_2, c_5$ are defined in Assumptions \ref{assumption::bounded_ra}-\ref{assumption::bounded_bw}, $K(t)$ is a time-varying gain matrix chosen as $K(t)=L_\gamma P(t)C(t)^\top Q(t)$, with $L_\gamma:=\textrm{blockdiag}(\gamma I_3,\gamma^2 I_3,\gamma^3 I_3)$ and $P(t)$ is solution to the following time-scaled continuous differential Riccati equation (CDRE):
\begin{align}
\label{eq:dP}
\frac{1}{\gamma}\dot P=AP+PA^\top -PC(t)^\top Q(t)C(t)P+V(t),
\end{align}
where $P(0)\in\mathbb{R}^{9\times 9}$ is positive definite, $Q(t)\in\mathbb{R}^{m\times m}$ and $V(t)\in\mathbb{R}^{9\times 9}$ are continuous, bounded and uniformly positive definite matrices. 
Furthermore, in the particular case where the output matrix $C(t)$ is constant ({\it e.g.,} range measurements as described in Section \ref{section:ranges}), it is possible to take $P$ as solution to the following continuous algebraic Riccati equation (CARE):
\begin{align}
AP+PA^\top -PC^\top QCP+V=0.
\end{align}
Selecting the gain $K$ offline might be very helpful to avoid the considerable computational burden associated with the real-time update of the CDRE.
\begin{figure}
    \centering
    \includegraphics[width=\columnwidth]{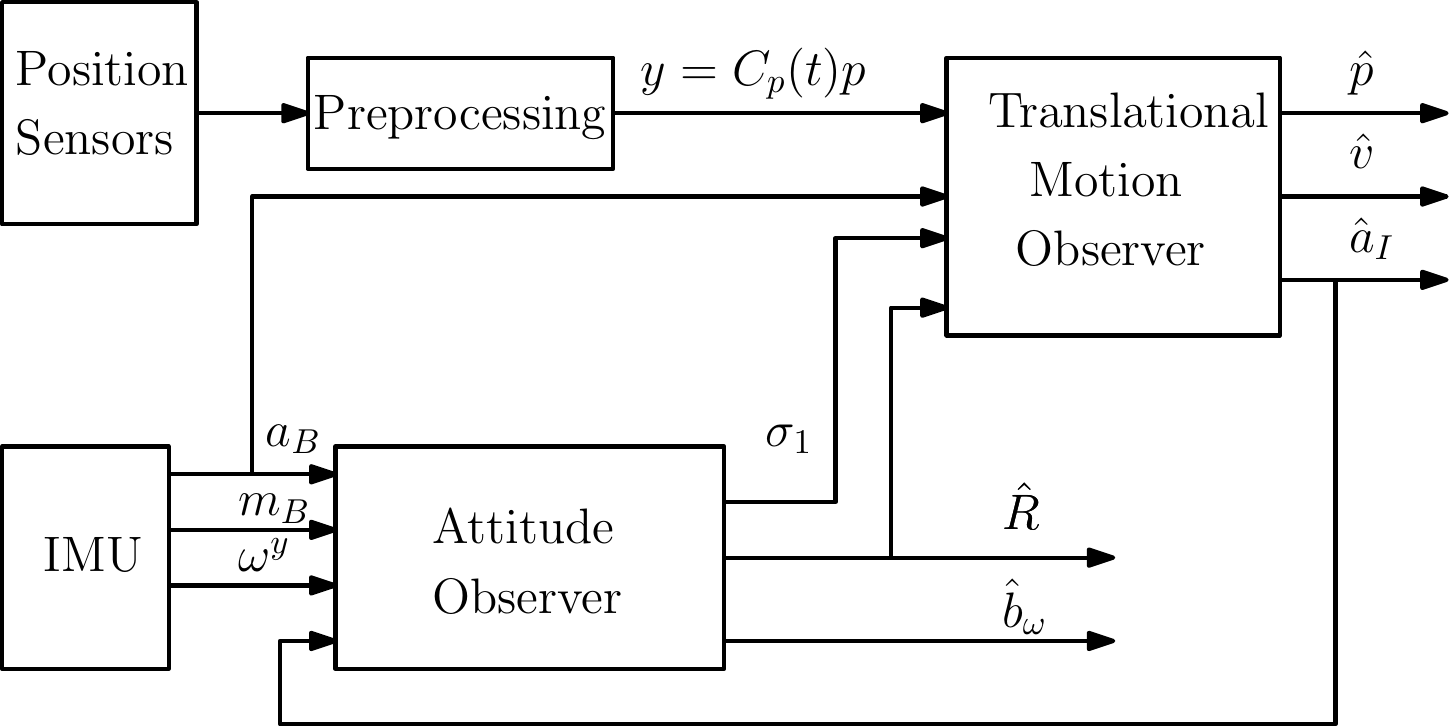}
    \caption{Overall structure of the proposed nonlinear observer. The observer uses an IMU and generic position sensors. The position information is first prepossessed to yield a linear output vector $y=C_p(t)x$.}
    \label{fig:structure}
\end{figure}
The overall strucutre of the proposed estimation algorithm is depicted in Fig.~\ref{fig:structure}. Before we state the stability result of the proposed observer, some remarks are in order.

The proposed nonlinear observer consists of an interconnection of the complementary filter-like attitude estimator \eqref{eq:dRhat}-\eqref{eq:dbhat}, which is inspired from \citep{Mahony2008NonlinearGroup,Grip2013}, with a high-gain observer \eqref{eq:dxhat}-\eqref{eq:dzhat} for the translational state estimation. Recall that the attitude complementary filter requires body-frame measurements of at least two non-collinear vectors known in the inertial frame. In applications involving low linear accelerations, one can directly use the accelerometer and magnetometer measurements in the filter. However, in the case of accelerated rigid body systems, this is not feasible as the acceleration in the inertial frame $a_I$ is unknown and cannot be approximated by the gravity vector. In this situation, we use an estimate of the inertial acceleration  $\hat a_I=B_2^\top\hat x$ obtained from the high-gain translational observer. The extra term $\sigma_1$, which vanishes as $\hat R$ tends to $R$, is used to yield a nice bounded coupling term between the tranlational and rotational error dynamics.
Furthermore, note that if we assume $\hat R\equiv R$ in \eqref{eq:dxhat}-\eqref{eq:dzhat},  then $\hat x$ satisfies $\dot{\hat x}=A\hat x+B_1ge_3+B_2\dot a_I+K(t)(y-C(t)\hat x)$ which is the Luenberger observer for the translational dynamics \eqref{eq:dx}-\eqref{eq:y} in the case where $\dot a_I$ is known. The introduction of the auxiliary state $\hat z$ allows to deal with the unknown input $B_2 \dot a_I$ of the translational dynamics \eqref{eq:dx} through the addition, after the integration of \eqref{eq:dzhat}, of the term $B_2\hat{R}a_B=B_2\tilde{R}^Ta_I$ which coincides with $B_2a_I$ when $\tilde{R}=I$.

The gain $K(t)$, inspired from \citep{Johansen2017NonlinearMeasurements}, relies on the high gain matrix $L_\gamma$ commonly used in high-gain observers (see, for instance, \citep{esfandiari1992output,saberi1990}). The high-gain observer for the translational dynamics allows to deal efficiently with the interconnection between the rotational and translational dynamics. In fact, as it is going to be shown later, this design will help to force the translational estimation errors to converge arbitrarily fast to a certain ball (which can be made arbitrary small) via a high gain parameter $\gamma$, preventing the attitude estimation error from converging to the undesired equilibria.  Consequently, asymptotic stability of the overall interconnected translational and rotational dynamics is guaranteed semi-globally. Finally, note that choosing $\gamma=1$ in \eqref{eq:dP} yields the traditional continuous differential Riccati equation (CDRE) used in the Kalman filter where $Q^{-1}(t)$ and $V(t)$ are interpreted as the covariance matrices for the output $y$ and the process, respectively. 
\subsection{Error Dynamics and Stability Analysis}\label{section:stability}
Let us define the following estimation error variables:
\begin{align}
\label{eq:tildex}
\tilde x&:=x-\hat x,\\
\tilde R&:=R\hat R^\top,\\
\label{eq:tildeb}
\tilde b_\omega&:=b_\omega-\hat b_\omega.
\end{align}
The geometric attitude estimation error $\tilde R$ defined above has been widely used in the literature of attitude observer design on $\mathbb{SO}(3)$; see, for instance, \citep{Mahony2008NonlinearGroup,Berkane2017OnSO3} among others. In view of \eqref{eq:dx}-\eqref{eq:y}, \eqref{eq:dxhat}-\eqref{eq:dRhat}, and using the fact that $AB_2=B_1$, the closed-loop system can be written as follows:
\begin{align}
				  \label{eq:dx_closed}
\dot{\tilde x}&=(A-K(t)C(t))\tilde x+B_2g(t,\tilde R,\tilde b_\omega),\\
\dot{\tilde R}&=\tilde R[-\hat R(\tilde b_\omega+k_1\sigma_2)]_\times,\\
\dot{\tilde b}_\omega&=\mathbf{P}_{c_5}^{\epsilon_b}(\hat b_\omega,k_2\sigma_2),
\end{align}
where we have defined $g(t,\tilde R,\tilde b_\omega):=(I-\tilde R)^\top\dot a_I(t)+\tilde R^\top [a_I(t)]_\times R(t)\tilde b_\omega$. Let us now introduce the time-scaled estimation error $\zeta:=L_\gamma^{-1}\tilde x$; see, for instance \citep{Johansen2017NonlinearMeasurements,esfandiari1992output}. In view of \eqref{eq:dx_closed}, the new error variable $\zeta$ satisfies the following time-scaled dynamics
\begin{equation}
\label{eq:dzeta}
\begin{aligned}
\frac{1}{\gamma}\dot\zeta
&=(A-P(t)C(t)^\top Q(t)C(t))\zeta+\frac{1}{\gamma^4}B_2g(t,\tilde R,\tilde b_\omega),
\end{aligned}
\end{equation}
where we have used the following easy-to-check facts:
\begin{align*}
L_\gamma^{-1}AL_\gamma=\gamma A,\quad L_\gamma^{-1}B_2=\gamma^{-3}B_2,\quad
C(t)L_\gamma=\gamma C(t).
\end{align*}
The first term of \eqref{eq:dzeta} corresponds to the nominal time-scaled dynamics of $\zeta$. It corresponds exactly to the closed-loop dynamics of a Bucy-Kalman estimator where $P(t)$ is designed via a differential Riccati equation. This motivates the time-scaling involved in the  Riccati equation in \eqref{eq:dP} to match the dynamics of $\zeta$ as done in \citep{Johansen2017NonlinearMeasurements}. To state the stability result of the observer, we further define the global estimation error variable 
$\varsigma:=[|\tilde R|,\|\tilde b_\omega\|,\|\zeta\|]^\top$.
\begin{thm}\label{theorem:navigation}
Consider the interconnection of the dynamics \eqref{eq:dp}-\eqref{eq:dR} with the observer \eqref{eq:dxhat}-\eqref{eq:sigma2:2} where Assumptions \ref{assumption::obsv}-\ref{assumption::bounded_bw} are satisfied. Assume, moreover, that the pair $(A,C(\cdot))$ is uniformly observable. Then, for each $\epsilon\in(\frac{1}{2},1)$ and $T>0$ and for all initial conditions such that $\zeta(0)\in\mathbb{R}^6$, $\tilde R(0)\in\mho(\epsilon)=\{\tilde R: |\tilde R(0)|\leq\epsilon\}$ and $\hat b_\omega(0)\in\mathbb{B}_{c_5+\epsilon_b}$, there exist $k_1^*>0$ and $\gamma^*\geq 1$ such that, for all $k_1\geq k_1^*$ and $\gamma\geq\gamma^*$, the estimation error $\varsigma(t)$ is globally uniformly bounded and
\begin{align}
\|\varsigma(t)\|\leq k\exp(-\lambda(t-T))\|\varsigma(T)\|\quad\forall t\geq T,
\end{align}
for some positive scalars $k$ and $\lambda$.
\end{thm}
\begin{pf}
See Appendix \ref{proof:theorem:navigation}.
\end{pf}
Theorem \ref{theorem:navigation} shows that the proposed nonlinear navigation observer guarantees exponential stability of the zero estimation error provided that the initial conditions of the estimation errors lie inside a compact set which can be arbitrary enlarged by an adequate tuning of the gains. Note that the gains conditions, provided in the proof, are rather conservative and simulation/experimental results have shown that the proposed navigation estimator has a large region of attraction regardless of the choice of the positive gains. It should be mentioned, that due to topological considerations, it is impossible to achieve global asymptotic (exponential) stability results with any continuous time-invariant attitude observer on $\mathbb{SO}(3)$.
Finally, the important condition on the uniform observability of the pair $(A,C(\cdot))$ will be discussed in the next section for different application scenarios.

\subsection{Observability Analysis}\label{section:observability}
A necessary condition for the result of Theorem \ref{theorem:navigation} to hold is the uniform observability of the pair $(A,C(\cdot))$. This condition intuitively means that the measurement of $y$ in \eqref{eq:yp} is enough to construct a converging translational observer (assuming perfect knowledge of the jerk $\dot a_I$) for \eqref{eq:dx}-\eqref{eq:y}. Now we derive the following important lemma. 
\begin{lemma}\label{lemma:obsv:general}
$(A,C(\cdot))$ is uniformly observable  if and only if there exist $\delta,\mu>0$ such that for all $t\geq 0$ one has
\begin{align}\label{condition:pe:general}
\frac{1}{\delta}\int_t^{t+\delta}C_p^\top(s)C_p(s)ds\geq\mu I_3.
\end{align}
\end{lemma}
\begin{pf}
See Appendix \ref{appendix:lemma:obsv:general}
\end{pf}
Lemma \ref{lemma:obsv:general} provides a persistency of excitation (PE) condition on the position output matrix $C_p(t)$ which is equivalent to the required uniform observability of the pair $(A,C(\cdot))$. In the case of constant $C_p$ the condition is equivalent to $\mathrm{rank}(C_p)=3$. This result has an intuitive meaning: uniform observability of the position subsystem (single integrator) is equivalent to uniform observability of the whole translational subsystem (multiple integrators). In the following subsections we will discuss this observability condition for different sets of positions sensors (range, bearing,...etc).

\subsubsection{Full Position Measurements}
Most GPS receivers provide position estimation when an unobstructed line of sight to four or more GPS satellites exists ({\it e.g.,} outdoor scenarios). In this case, the full position is assumed available and the output equation \eqref{eq:yp} is taken with $C_p(t)=I_3$.  In this case, the pair $(A,C)$ is uniformly observable since $\textrm{rank}(C_p)=3$. Different algorithms such as \citep{Grip2013,Bryne2017NonlinearAspects} have developed nonlinear observers relying on full position measurements (loosely coupled integration). Note that altitude determined using low-cost GPS is not generally reliable enough which can motivate the use of a pressure altimeter to determine the altitude.

\subsubsection{Range Measurements}\label{section:ranges}
Different sensors can be used to provide range measurements. For instance, in a tightly coupled GPS/INS integration, raw GPS observations (range measurements) are used directly in the estimation scheme to allow the use of fewer observations than actually needed to reconstruct the position. Another instance of range sensors is the UWB technology which has proved successful especially in indoor applications. 

Assume that we have available $n$ source points (anchors) with known possibly time-varying\footnote{Though in most practical applications the anchors for range measurements are at constant locations.} locations (positions), denoted as $p_i$. The corresponding range measurements are given by
\begin{align}\label{eq:di}
d_i=\|p-p_i\|,\quad i=1,\cdots,n.
\end{align}
To obtain an output equation of type \eqref{eq:yp} we proceed as follows, see \citep{Hamel2017PositionMeasurements}. Let $\bar y_i:=0.5(d_i^2-\|p_i\|^2)$ and define the weighted output
\begin{equation}
    \bar y_0:=\sum_{i=1}^n\alpha_i\bar y_i,
\end{equation}
where $\alpha=[\alpha_1,\cdots,\alpha_n]\in\mathbb{R}^n$ is a vector of constant real numbers such that $\textstyle\sum_{i=1}^n\alpha_i=1$. Now, we define our output vector as follows:
\begin{align}\label{eq:Cp:range}
    y:=\begin{bmatrix}
    \bar y_1-\bar y_0\\
    \vdots\\
    \bar y_n-\bar y_0
    \end{bmatrix}=\begin{bmatrix}
    \bar p_1^\top\\
    \vdots\\
    \bar p_n^\top
    \end{bmatrix}p:=C_p(t)p
\end{align}
where we have defined the following known vectors
\begin{align}\label{eq:barpi}
    \bar p_j:=\sum_{i=1}^n\alpha_i(p_i-p_j),\quad j=1,\cdots,n.
\end{align}
Output equation \eqref{eq:Cp:range} encodes the position information provided by the range measurements. If we have an additional altimeter sensor, the output of the altimeter can be written as $h=e_3^\top p$. In this case, and to include this measurement, we can add the row vector $e_3^\top$ to the matrix $C_p(t)$ defined in \eqref{eq:Cp:range}.The following lemma then immediately follows from Lemma \ref{lemma:obsv:general} and \eqref{eq:Cp:range}.
\begin{lemma}\label{lemma:obs:range}
$(A,C(\cdot))$ is uniformly observable if and only if there exist $\delta,\mu>0$ such that for all $t\geq 0$ one has
\begin{align}\label{condition:pe}
\frac{1}{\delta}\int_t^{t+\delta}\sum_{i=1}^n\bar p_i(s)\bar p_i^\top(s)ds+\alpha e_3e_3^\top\geq\mu I_3.
\end{align}
where $\alpha=1$ if the altimeter is used or $\alpha=0$ otherwise.
\end{lemma}

Different conclusions can be drawn from the observability condition in Lemma \ref{lemma:obs:range}. In the case where $n=1$, one has $\bar p_1=0$ and therefore the condition is not fulfilled. When $n=2$ and $\alpha=0$,  condition \eqref{condition:pe} is equivalent to a P.E. condition on the vector $p_1(t)-p_2(t)$. Roughly speaking, the latter condition prevents the vector $p_1(t)-p_2(t)$ from staying indefinitely in any plane. If we add an altimeter, \textit{i.e.,} $n=2$ and $\alpha=1$, the condition prevents the vector $p_1(t)-p_2(t)$ from staying indefinitely in the plane containing $e_3$. When $n=3$ and $\alpha=0$, condition \eqref{condition:pe} is not satisfied if the three anchors $p_1, p_2$ and $p_3$ are in the same plane for all times. However, when adding an altimeter, the condition is satisfied when the three anchors $p_1, p_2$ and $p_3$ are not aligned and $e_3$ does not belong to the plane spanned by the three anchors. Finally, when $n\geq4$ the condition holds if at least $4$ anchors are non-coplanar. The latter result is consistent with the minimum number of anchors required to geometrically (multilateration) find the position with no ambiguity from multiple range measurements (for stationary anchors). Adding an altimeter sensor relaxes the requirement  of $4$ non-coplanar anchors to $3$ non-aligned anchors. Note that in the case of constant anchors positions, the only two cases where the observability is guaranteed are: 1) the case of $n\geq 4$, and 2) the case of $n=3$ and $\alpha=1$.


Besides, we would like to emphasize that the observability condition in Lemma \ref{lemma:obs:range} is independent from the trajectory of the vehicle and depends only on the position of the anchors. In \citep{Hamel2017PositionMeasurements} the variable $\bar y_0$ was added, via state augmentation, to the overall state of the system and the derived observability condition was explicitly written as a PE condition on the input. This allows to design the observer even when $n=1$ (single range measurement) under a PE condition on the input. However the focus of \citep{Hamel2017PositionMeasurements} was on position estimation only, whereas here we consider the full navigation problem. Handling the case of single range measurements, via state augmentation, for the full navigation problem exceeds the scope of this paper.

\subsubsection{Bearing Measurements}\label{section:bearings}
A motion capture system consisting of a stationary array of cameras capturing the vehicle from multiple angles can be used to provide raw bearing measurements in the inertial frame. Each camera can provide a unit vector direction measurement as follows:
\begin{align}\label{eq:bearings}
y_i=R_i^\top\frac{p-p_i}{\|p-p_i\|},\quad i=1,\cdots,n,
\end{align}
where $R_i\in\mathbb{SO}(3)$ is the known orientation of the $i$-th camera with respect to the inertial frame, $p_i$ is the inertial position of the $i$-th camera and $n$ is the number of working cameras. Note that here we consider bearings measured in the inertial frame in contrast to \citep{Hamel2018RiccatiProblem} for examples where bearings are measured in the body frame of reference. Inspired from \citep{Batista2015NavigationMeasurements,Hamel2017PositionMeasurements}, we generate an artificial output, which is linear in the position, by noticing that
	\begin{align}
		\Pi(y_i)R_i^\top(p-p_i)=0,
	\end{align}
where $\Pi:\mathbb{S}^2\to\mathbb{R}^{3\times 3}$ is the orthogonal projection map defined as
\begin{align}\label{eq:pi}
\Pi(z):=I_3-zz^\top.
\end{align}
Consequently,  it is possible to define the following output vector:
\begin{align}\label{eq:Cp:bearings}
y=\begin{bmatrix}
\Pi(y_1)R_1^\top p_1\\
\vdots\\
\Pi(y_n)R_n^\top p_n\\
\end{bmatrix}=\begin{bmatrix}
\Pi(y_1)R_1^\top\\
\vdots\\
\Pi(y_n)R_n^\top
\end{bmatrix}p:=C_p(t)p.
\end{align} 
If we have an additional altimeter sensor, the output of the altimeter can be written as $h=e_3^\top p$. In this case, and to include this measurement, we can add the row vector $e_3^\top$ to the $C_p(t)$ matrix defined in \eqref{eq:Cp:bearings}. 
In the case of the bearing measurements discussed above and a possible use of an altimeter, we state the following lemma which is an immediate consequence of Lemma \ref{lemma:obsv:general} and \eqref{eq:Cp:bearings}.
\begin{lemma}
$(A,C(\cdot))$ is uniformly observable if and only if there exist $\delta,\mu>0$ such that for all $t\geq 0$ one has
\begin{align}\label{condition:pe:bearings}
\frac{1}{\delta}\int_t^{t+\delta}\sum_{i=1}^n\Pi(R_iy_i(s))ds+\alpha e_3e_3^\top\geq\mu I_3.
\end{align}
where $\alpha=1$ if the altimeter is used or $\alpha=0$ otherwise.
\end{lemma}%
Depending on the number of bearing measurements available and whether or not an altimeter sensor is used, we derive the following cases that fulfill the required observability condition.
\begin{lemma}\label{lemma:obsv}
Assume that the velocity $v$ is bounded. Condition \eqref{condition:pe:bearings} is satisfied if one of the following conditions holds:
\begin{itemize}
\item[i)] There exist $i,j \in \{1,\hdots,n\}$ and $\epsilon>0$, such that for all $t^*>0$, there exists $t>t^*$, such that $\|R_iy_i(t)\times R_jy_j(t^*)\|\geq\epsilon$.
\item [ii)] at least three cameras are not aligned.
\end{itemize}
Moreover, if the altimeter is used, condition \eqref{condition:pe:bearings} holds if either (i) or (ii) are satisfied, or one of the following conditions holds:
\begin{itemize}
    \item[iii)] There exist $i\in \{1,\hdots,n\}$, $\epsilon>0$, such that for all $t^*>0$, there exists $t>t^*$, such that $|e_3^\top R_i y_i(t)|\geq\epsilon$.
    \item[iv)] There exist $i, j \in \{1,\hdots,n\}$, such that $e_3^\top(p_i-p_j)\neq 0$.
\end{itemize}
\end{lemma}
\begin{pf}
See Appendix \ref{appendix:lemma:obsv}.
\end{pf}
Different conclusions can be derived from Lemma \ref{lemma:obsv}. First, regardless whether we use an altimeter or not, if $v$ is bounded, uniform observability is satisfied if we use
\begin{itemize}
    \item at least one bearing measurement $y_i$ which is not constant for all times. In other words, this requires that the vehicle is never static nor indefinitely moving in a straight line with $p_i$.
    \item at least two bearing measurements $y_i$ and $y_j$ which are not aligned for all times. This means that the vehicle is not indefinitely aligned with both cameras.
    \item at least three non-aligned cameras, regardless of the trajectory of the vehicle.
\end{itemize}
In addition, if we use the altimeter measurement, uniform observability is satisfied if
\begin{itemize}
    \item at least one bearing measurement $y_i$ is available, and the vehicle is not indefinitely located at the same altitude as $p_i$.
    \item  at least two cameras are not at the same altitude, regardless of the vehicle's trajectory.
\end{itemize}
Interestingly, the use of an altimeter allows to guarantee uniform observability for all trajectories with only two cameras installed at different altitudes. This fact can be explained as follows: if the two cameras are at the same altitude and the vehicle is moving on a trajectory passing through both cameras, then all the three measurements (two bearings and altimeter) will be constant. Consequently, it is not possible to observe the position of the vehicle since its trajectory does not affect the sensor readings. On the other hand, if the two cameras are not at the same altitude, the previous scenario is not possible since any trajectory which is aligned with both cameras would necessarily cause the altitude to change. Finally, for planar trajectories at a constant altitude, uniform observability is guaranteed with the use of one camera located at a different altitude than the altitude of the flight.

\section{Simulation Results}\label{section:simulation}
In this section, we simulate the nonlinear observer of Section \ref{section:riccati_observer} with different position sensing scenarios. In all simulations, the angular velocity applied to the vehicle is given by:
\begin{align}
\omega(t)=\begin{bmatrix}
\sin(0.1t+\pi)\\
0.5\sin(0.2t)\\
0.1\sin(0.3t+\pi/3)
\end{bmatrix}
\end{align}
with an initial attitude $R(0)=\exp([\pi e_2]_\times/2)$.  The gyro measurements are corrupted by a constant bias of $2\;\mathrm{(deg/sec)}$ in each axis. The inertial earth's magnetic field is taken as $m_I=[0.033\;0.1\;0.49]^\top$ and the earth's gravity is $g=9.81\;(\mathrm{m/sec}^2)$.

\subsection{Range Measurements}
\begin{figure}
    \centering
    \includegraphics[clip,trim = 3mm 12mm 0mm 03mm, width=.44\textwidth]{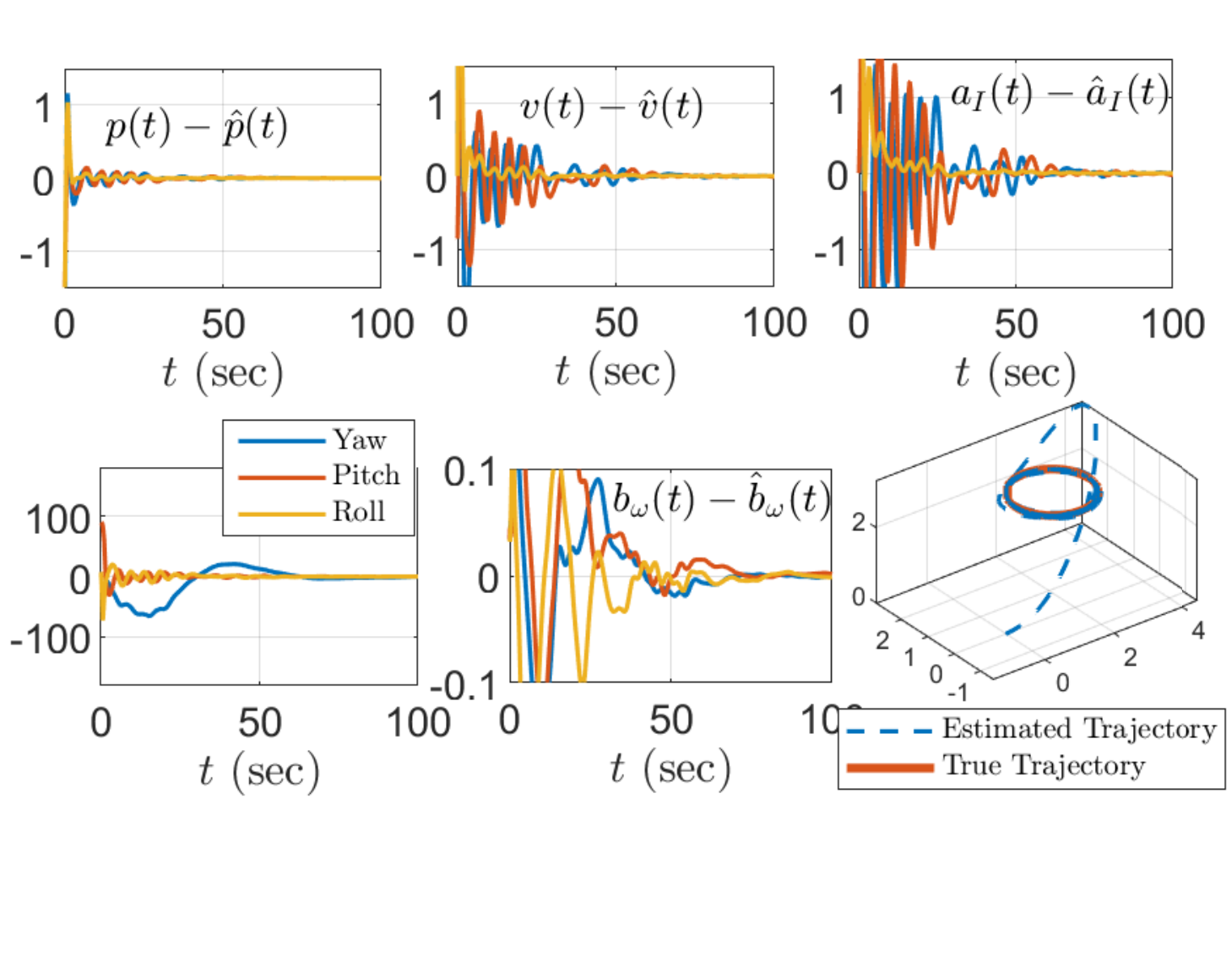}
    \caption{Estimation errors and trajectory in the case of 4 non-coplanar range measurements.}
    \label{fig:estimation:range}
\end{figure}

Here we assume available four non-coplanar source points $p_1,\cdots,p_4$ located at $p_1=0$ and $p_{i+1}=e_i, i=1,2,3$. We consider a vehicle moving along the circular trajectory 
\begin{align}
p(t)=\begin{bmatrix}
1.075\cos(\pi t/4)+2.5\\
1.075\sin(\pi t/4)+1.5\\
2.2
\end{bmatrix}.
\end{align}
 The initial conditions for the observer states are $\hat z(0)=\hat b_\omega(0)=0, \hat R(0)=I_3,$ and $P(0)=I_9$. The parameters of the observer are selected as $k_1=2, k_2=1, \rho_1=1, \rho_2=0.1,\epsilon_b=0.001, c_5=0.06, \hat c_2=15, \gamma=2, V(t)=I_3$ and $Q=5I_4$. The simulation results given in Figure \ref{fig:estimation:range} show that the proposed observer was able to estimate the position, velocity, acceleration, attitude and gyro bias using IMU and range measurements.

\subsection{Bearing Measurements}
\begin{figure}
\centering
\includegraphics[clip,trim = 7.5mm 20mm 15mm 7mm, width=.47\textwidth]{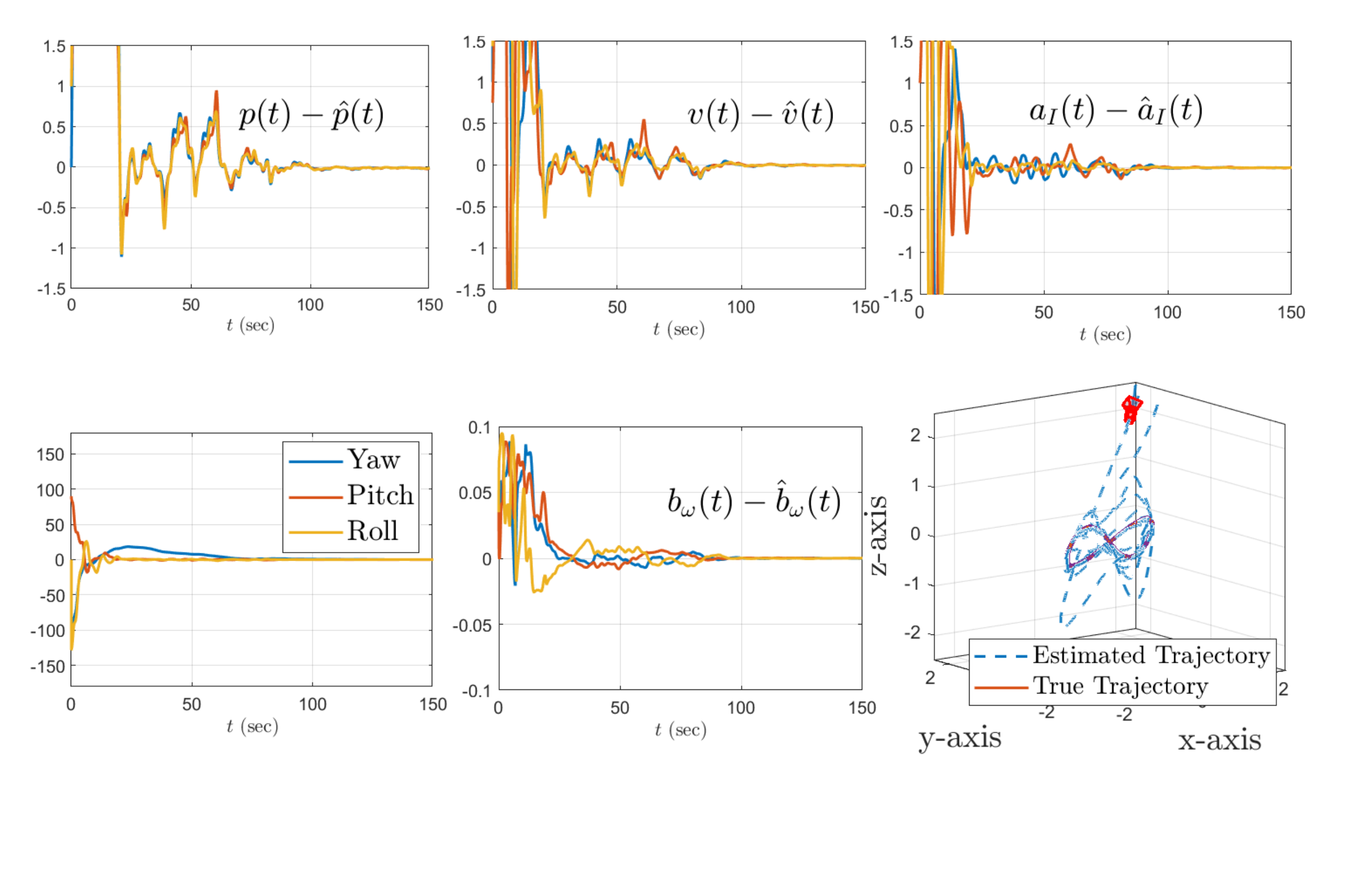}
\caption{Estimation errors and trajectory in the case of a single bearing measurement. }
\label{figure1}
\end{figure}

\begin{figure}
\centering
\includegraphics[clip,trim = 7.5mm 20mm 15mm 7mm, width=.47\textwidth]{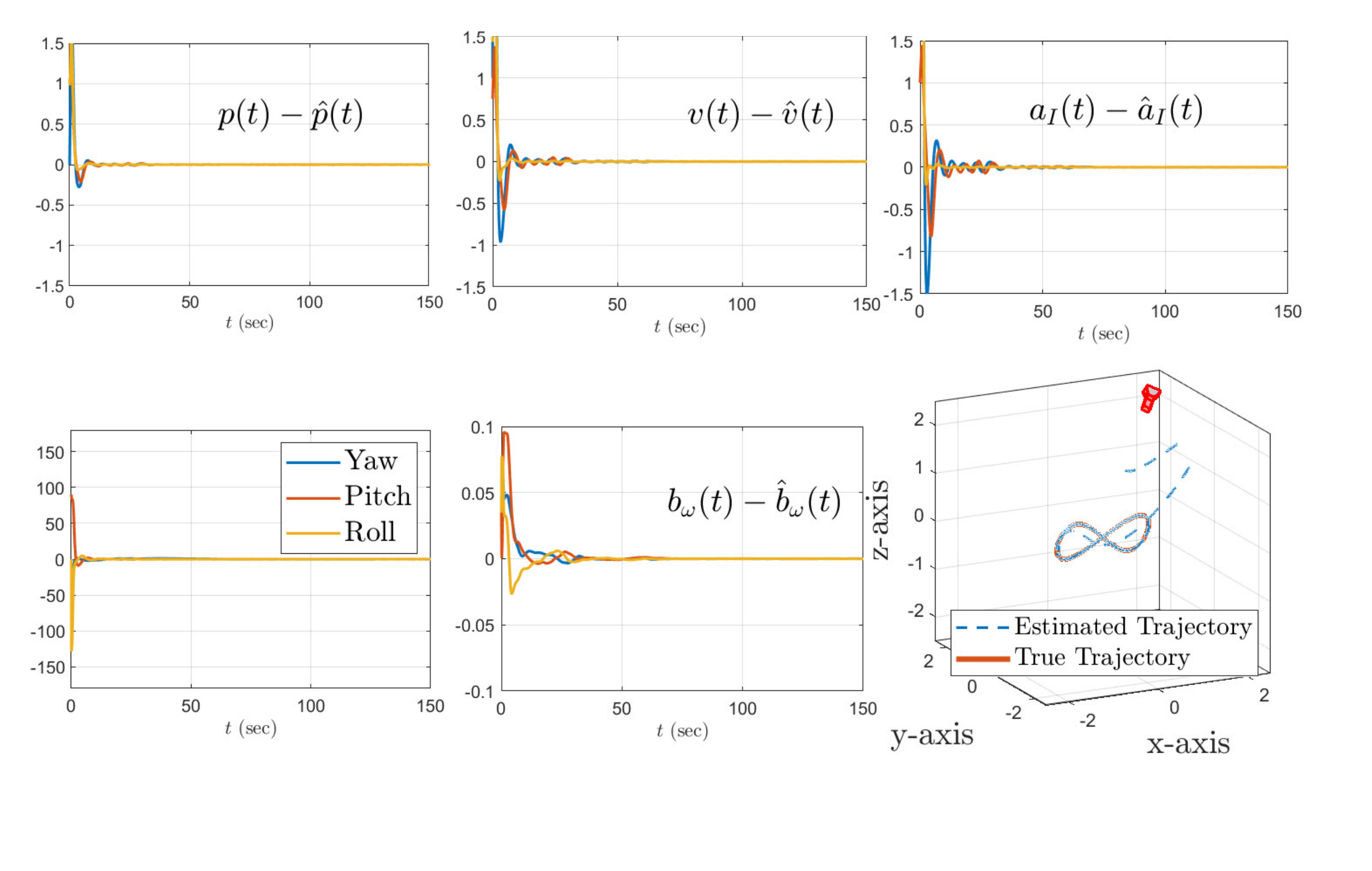}
\caption{Estimation errors and trajectory in the case of a single bearing measurement and an altimeter sensor. }
\label{figure2}
\end{figure}

Here we consider inertial frame bearing measurements as described in Subsection \ref{section:bearings}. Consider a vehicle moving on following eight-shaped trajectory:
\begin{align}
p(t)=\begin{bmatrix}
\cos(t/2)\\
\sin(t)/4\\
-\sqrt{3}\sin(t)/4
\end{bmatrix}.
\end{align}
The initial conditions for the observer states are $\hat z(0)=[1\;1\;1\;1\;1\;1\;1\;1\;1]^\top, \hat R(0)=I_3, \hat b_\omega(0)=[0\;0\;0]^\top$ and $P(0)=I_9$. The parameters of the observer are selected as $k_1=2, k_2=1, \rho_1=1, \rho_2=0.1,\epsilon_b=0.001, c_5=0.06, \hat c_2=15, \gamma=2, V(t)=I_3$ and $Q=I_m$. We consider two simulation scenarios. In the first scenario, we consider a single camera, located at $p_1=[2,2,2]^\top$ pointing towards the origin, providing a single bearing measurement. In the second scenario, we consider a single bearing measurement along with an altimeter.  Figures \ref{figure1} and \ref{figure2} show the evolution of the different estimation errors versus time. In both scenarios, the estimation errors converge to zero. However, adding the altimeter sensor in the second scenario has considerably improved the convergence rate compared to the first simulation scenario without an altimeter.

%
\section{Experimental Results}\label{section:experiments}
In this section, we experimentally validate the proposed nonlinear navigation observer on a dataset recorded with a custom Quadrotor platform.

\subsection{Experimental Setup}
\begin{figure}
    \centering
    \includegraphics[clip,trim = 00mm 0mm 00mm 00mm, width=.44\textwidth]{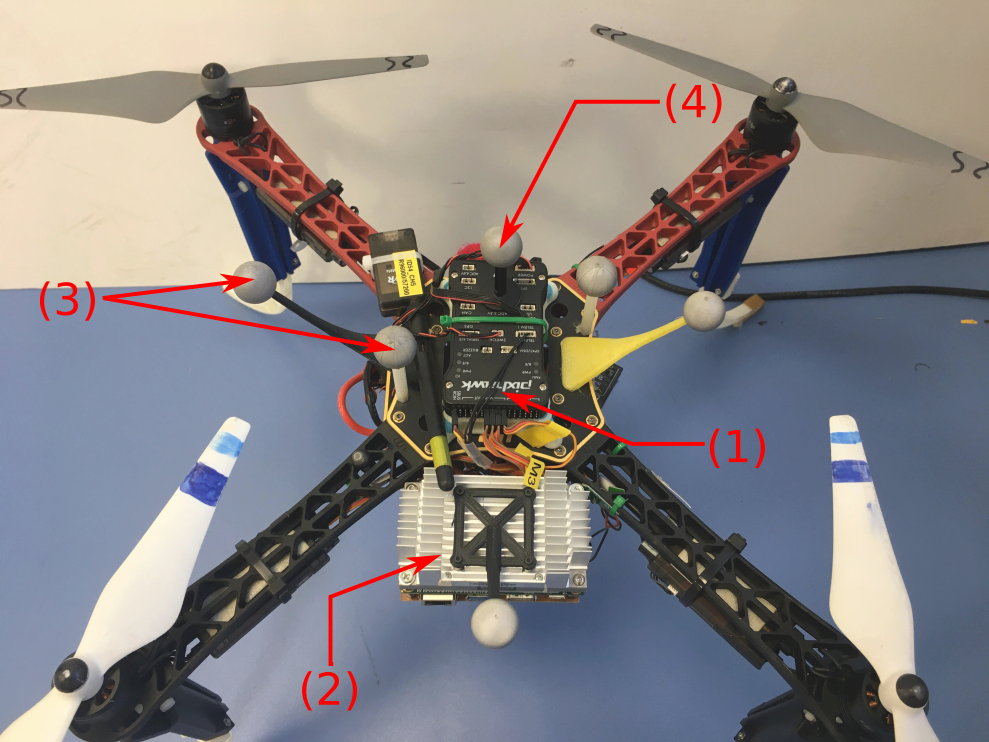}
    \caption{Quadrotor vehicle used in the experiments. (1) Pixhawk Flight Controller Unit (FCU). (2) Companion computer Nvidia Jetson TX1. (3) Optical Markers for the motion capture system (MOCAP). (4) Optical Marker used for bearings computations.}
    \label{figure:Quadrotor Experimental}
\end{figure}
The vehicle (see Figure \ref{figure:Quadrotor Experimental}), based on a DJI-450 frame, is equipped with a Pixhawk flight controller unit (PX4 software) and a Nvidia Jetson TX1 as companion computer. The Pixhawk unit is equipped with a 16 bit gyroscope (L3GD20H) and 14 bit accelerometer and magnetometer sensor (LSM303D), all the sensors measurement are sent to the companion computer using mavlink/mavros protocol. A custom trajectory tracking non-linear controller based on \citep{KAI2017} has been implemented on PX4 open source flight control software. 
\begin{figure}
    \centering
    \includegraphics[clip,trim = 08mm 05mm 7mm 05mm, width=.45\textwidth]{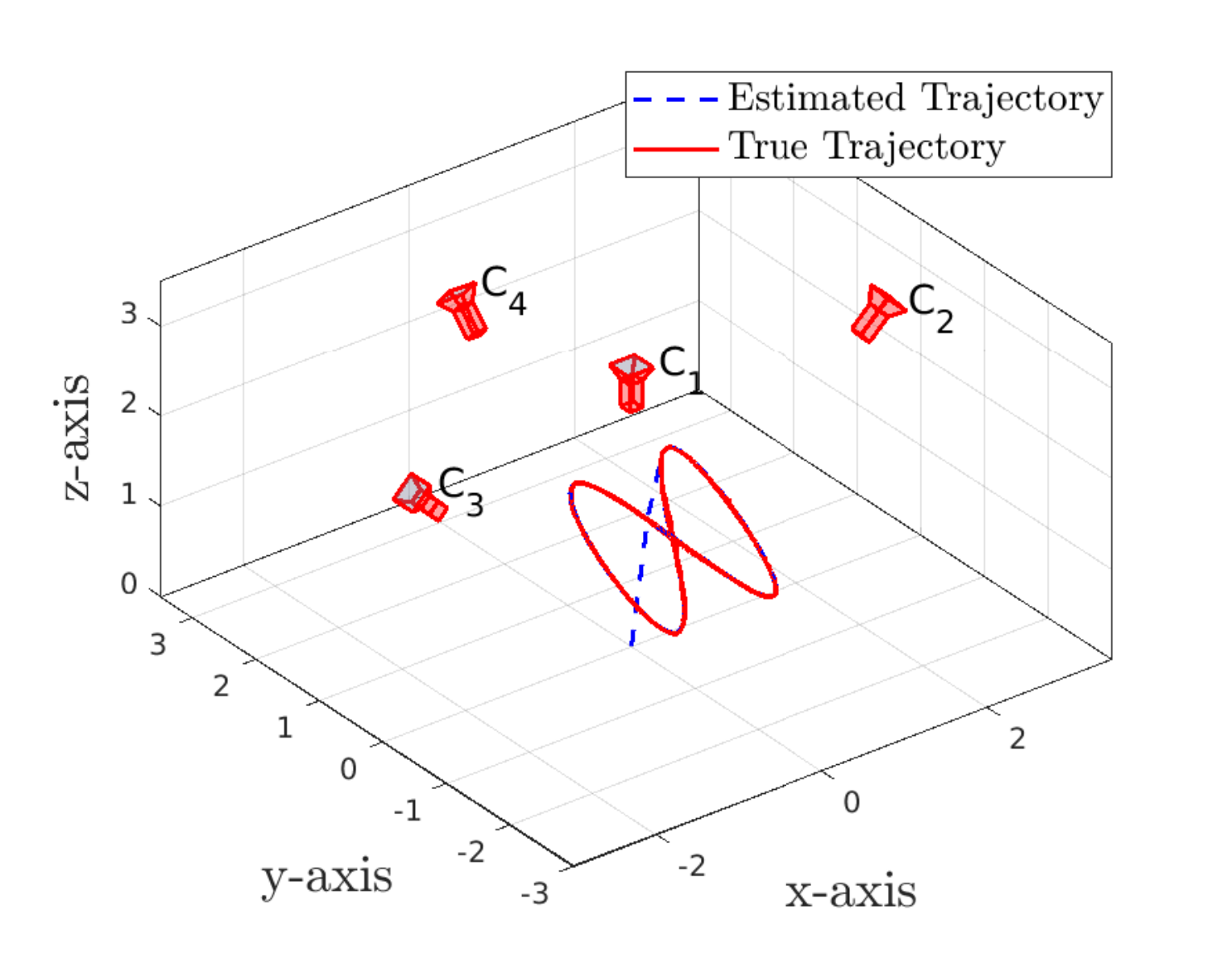}
    \caption{Position of the four cameras used for bearing measurements and 3D true trajectory vs estimated trajectory.}
    \label{figure:Experimental 3D trajectory}
\end{figure}

An OptiTrack motion capture system, comprising 8 cameras, is used along with optical markers mounted on the quadrotor (see Figure \ref{figure:Quadrotor Experimental}) in order to provide full pose ground truth measurements.
The linear velocity is obtained from the position via a filtered numerical derivative.
OptiTrack Camera SDK have been used to retrieve bearing measurements of the optical marker located on the top of the Pixhawk (Figure \ref{figure:Quadrotor Experimental} item (4)) from each camera.
A ground computer, connected to the motion capture system, sends ground truth measurement and bearings to the Nvidia TX1 over WiFi.

\subsection{Experiment and Results}
In the experiment, the quadrotor is commanded to track the following Lemniscate trajectory (depicted in Figure \ref{figure:Experimental 3D trajectory})
\[
    p_{ref}(t) = \begin{bmatrix} 
            0.8 \cos(w_{r} t)\\
            0.5+0.8 \sin(2 w_{r} t)\\
            1+0.8 \sin(2 {w_r} t)
            \end{bmatrix}
\]
with $w_{r}=0.16 \pi $, while maintaining a constant yaw angle of $-90$ degree. Notice that, due to the under-actuated nature of the vehicle, the roll and pitch angles can not be arbitrarily chosen and their time behaviours depend on the position control loop.

During the experiment, only four cameras have been used for bearing measurements. The inertial positions and orientations with respect to the inertial frame of the cameras (see Figure \ref{figure:Experimental 3D trajectory}) are the following:
\[
    \begin{array}{lllllllll}
        p_1 \!= \bigg[\begin{smallmatrix} 0\\0\\2.8 \end{smallmatrix}\bigg],&&  
        \!\!p_2 \!= \bigg[\begin{smallmatrix} 2.89\\0\\2.57 \end{smallmatrix}\bigg],&&
        \!\!p_3 \!= \bigg[\begin{smallmatrix} -2.44\\0\\2.42 \end{smallmatrix}\bigg],&& 
        \!\!p_4 \!= \bigg[\begin{smallmatrix} 0.08\\2.65\\2.37 \end{smallmatrix}\bigg], \\[1.5em]
        q_1\!=\bigg[\begin{smallmatrix} 0\\1\\0\\0 \end{smallmatrix}\bigg], && \!\!q_2\!=\bigg[\begin{smallmatrix} 0.5\\0\\-0.86\\0 \end{smallmatrix}\bigg], && \!\!q_3\!=\bigg[\begin{smallmatrix} 0.45\\0\\0.89\\0 \end{smallmatrix}\bigg], && \!\!\!\!\!q_4\!=\bigg[\begin{smallmatrix} 0.32\\0.63\\-0.63\\0.32 \end{smallmatrix}\bigg].
    \end{array}
\]

The parameters of the observer are selected as $k_1=10$, $k_2=1$, $\rho_1=4$, $\rho_2=1.0$, $\epsilon_b=0.001$, $c_5=0.3$, $\hat{c}_2=30$, $\gamma=2$, $V(t)=12 I_9$ and $Q=2 I_{12}$.
The initial conditions for the observer states are $\hat z(0)=0_{9\times1}$, $\hat{R}(0)=I_3$, $\hat{b}_\omega(0)=0_{3\times1}$ and $P(0)=10I_9$.
\begin{figure}
    \centering
    \includegraphics[clip,trim = 10mm 01mm 12mm 05mm, width=.46\textwidth]{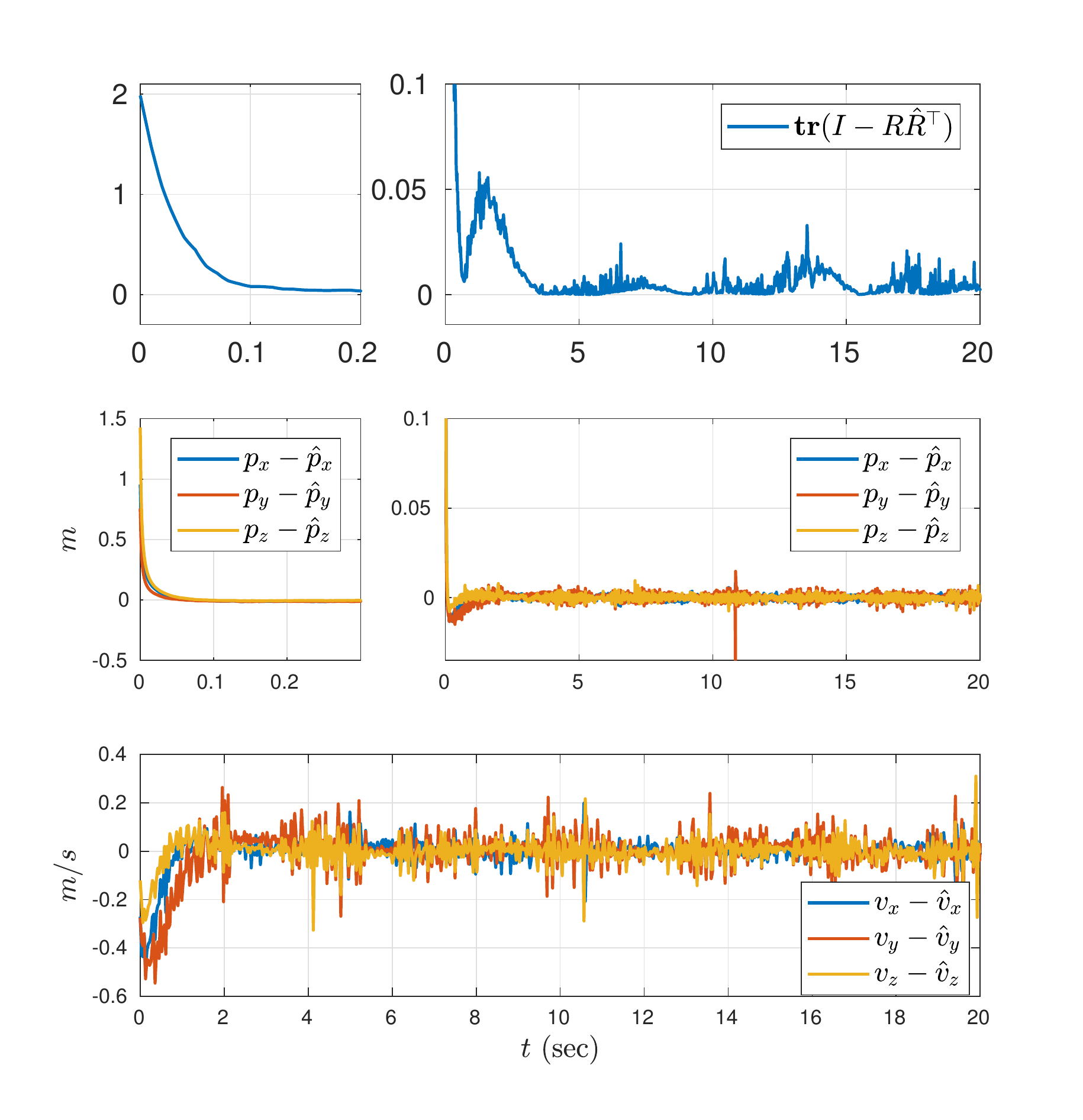}
    \caption{Estimation errors in case of 4 bearing measurements for the Quadrotor vehicle. }
    \label{figure:Experimental Error Plots}
\end{figure}

\begin{figure}
    \centering
    \includegraphics[clip,trim = 12mm 01mm 12mm 05mm, width=.44\textwidth]{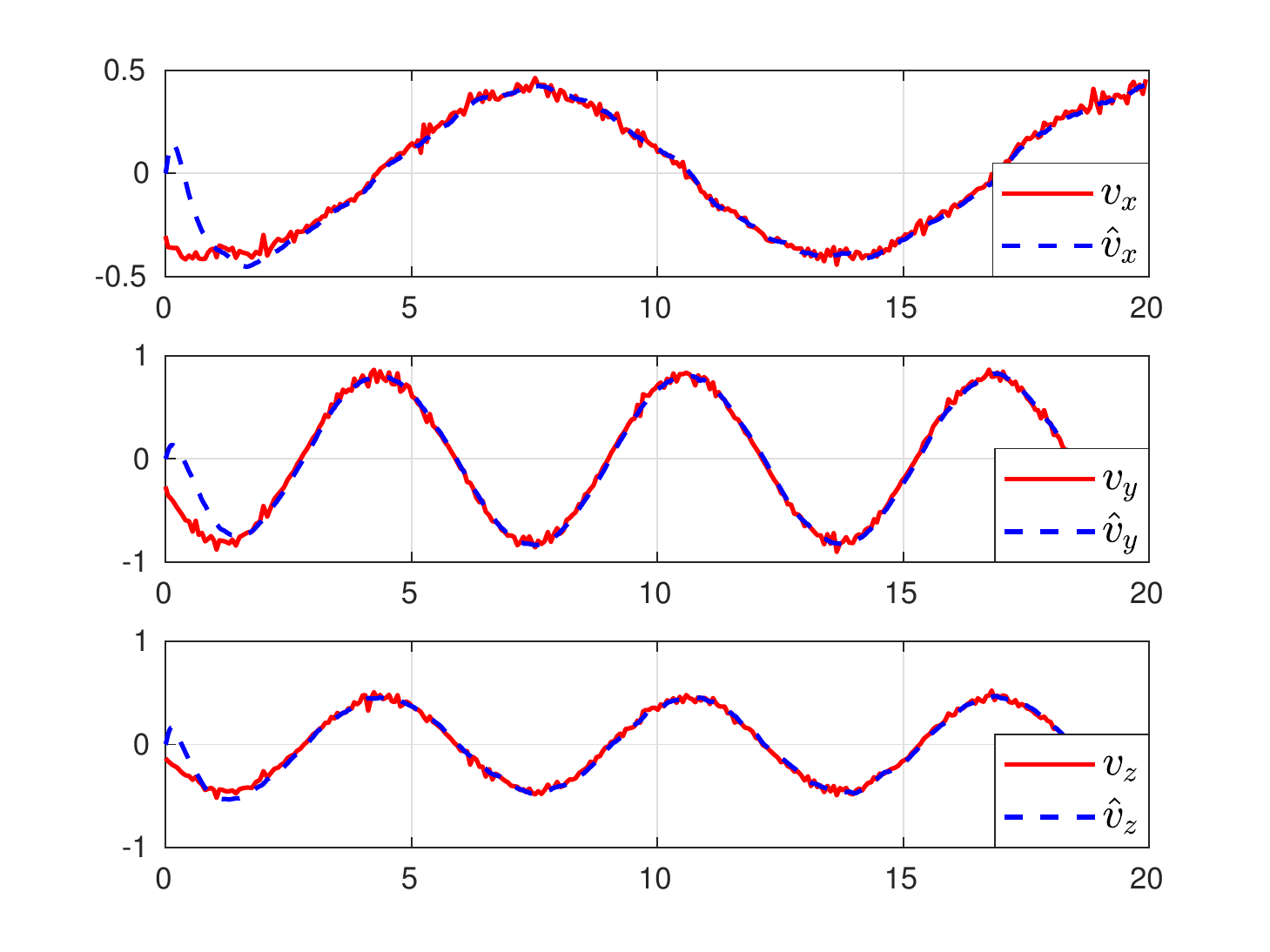}
    \caption{Time behaviour of the components of the real and estimated linear velocity.}
    \label{figure:Experimental true velocity vs estimated velocity}
\end{figure}
Figure \ref{figure:Experimental 3D trajectory} shows the 3D trajectory of the vehicle versus the observer estimated trajectory, whereas  Figure \ref{figure:Experimental Error Plots} shows the time evolution of the different estimation errors.
Plots clearly show that the estimated attitude $\hat R$, estimated position $\hat p$ and estimated velocity $\hat v$ converge to the real attitude $R$, real position $p$ and ground truth velocity $v$ of the quadrotor, respectively.Note that the noise in the velocity error (Figure \ref{figure:Experimental Error Plots}) is mainly due to ground truth obtained by numerical derivation as shown in Figure \ref{figure:Experimental true velocity vs estimated velocity}.

\subsection{Implementation Aspects}
Although the attitude observer in \eqref{eq:dRhat}-\eqref{eq:dbhat} is designed on $\mathbb{SO}(3)$, it is possible to lift its dynamics to the unit-quaternion space $\mathbb{S}^3$ as follows:
\begin{align}\label{eq:quaternionRepresentation}
\dot{\hat q}&= \frac{1}{2}\mathbf{T}(\hat\omega)~\hat q,
\end{align}
with $\hat\omega:=\omega^y-\hat b_\omega+k_1\sigma_2$ and
$$\mathbf{T}(\hat\omega):=\begin{bmatrix} 0 & -\hat{\omega}^\top \\ \hat{\omega} & -[\hat{\omega}]_\times \end{bmatrix}. $$
The use of the unit-quaternion representation (4 state components) for the observer implementation is more computationally efficient than the use of the rotation matrix representation (9 state components).


In our implementation, the attitude estimation kinematics in \eqref{eq:quaternionRepresentation} have been discretized as follows (\cite{whitmore2000closed}):
\begin{align*}
    \hat{q}_{k+1}&=\left[\cos\left(0.5 \| \hat{\omega}_k \| \tau \right)I_4+\frac{ \tau}{2} \mbox{sinc} \left(0.5 \| \hat{\omega}_k \| \tau \right) \mathbf{T}(\hat\omega) \right]\hat{q}_k,
\end{align*}
with
\[
\mbox{sinc}(\alpha)= \left\{ \begin{array}{ll}
1, & ~~\mbox{if}~~\alpha=0,\\
\frac{\sin \alpha}{\alpha}, & ~~\mbox{elsewhere}.
\end{array}
\right.
\]
and the sampling period $\tau=5ms$.
Note that this discretization  preserves the unit quaternion constraint and thus no normalization is required at each integration step.
The estimated rotation matrix can be retrieved from the unit-quaternion as $\hat{R}_k=\mathcal{R}(\hat{q}_k)$, where the Rodrigues map is given by $\mathcal{R}: \mathbb{S}^3 \mapsto \mathbb{SO}(3)$ such that $q:=[q_0,q_v^\top]^\top \rightarrow \mathcal{R}(q)=(q_0^2-||q_v||^2) I_3+2q_0[q_v]_\times+2q_v q_v^\top$. The translational estimation dynamics \eqref{eq:dxhat}-\eqref{eq:dzhat}  have been discretized using forward Euler integration with the same frequency of the IMU (200Hz). 

\section{Conclusion}\label{section:conclusion}
In this work, we proposed a nonlinear observer relying on IMU (accelerometer, gyroscope and magnetometer) measurements, and full/partial position information, for the simultaneous estimation of the position, velocity, orientation, and gyro bias. It employs a Riccati-like gain update which allows to possibly include the noise characteristics. The orientation estimates are obtained directly on the Special Orthogonal group of rotations.  The stability of the closed-loop system is proved to be exponential with a domain of attraction that can be arbitrary enlarged. A detailed observability analysis has been carried out for different types of measurements (full position, ranges, bearings, altimeter), and conditions on the number and location of the sensors with respect to the vehicle, guaranteeing the viability of the proposed observer, have been derived.  It is shown that the use of an altimeter relaxes the observability requirements, and enhances the observer's convergence rates. For example, in the case of two bearings and an altimeter, uniform observability is guaranteed for all trajectories as long as the two cameras are installed at different altitudes.

\appendix
\section{Proof of Theorem \ref{theorem:navigation}}\label{proof:theorem:navigation}
Thanks to assumption that $(A,C(\cdot))$ is uniformly observable and the fact that $Q(t)$ and $V(t)$ are uniformly positive definite and bounded matrices, there exist (see \citep{Bucy1972TheBounds}) positive scalars  $\beta_1$ and $\beta_2$ such that $\beta_1 I\leq P(t)\leq\beta_2 I$ for all $t\geq 0$. Note that these scalars are independent of $\gamma$, see \cite[Lemma 6]{johansen2015nonlinear}. Furthermore, in view of Assumption \ref{assumption::bounded_bw} and property \eqref{ineq:P1} of the projection mechanism, the bias estimation error is bounded such that $\|\tilde b_\omega\|\leq 2c_5+\epsilon_b:=c_b$. It follows, using Assumption \ref{assumption::bounded_ra}, that
\begin{equation}
\begin{aligned}
\|g(t,\tilde R,\tilde b_\omega)\|&\leq\sqrt{8}c_3|\tilde R|+c_2\|\tilde b_\omega\|\\
											   &\leq\sqrt{8}c_3+c_2c_b:=c_g.
											   \label{ineq:g}
\end{aligned}  
\end{equation}
Furthermore, in view of \eqref{eq:dP}, the time derivative of $P^{-1}(t)$ satisfies
\begin{equation}
\label{eq:dP-1}
\begin{aligned}
&\frac{1}{\gamma}\dot{P}^{-1}\\
&=-P^{-1}A-A^\top P^{-1}+C^\top(t)Q(t)C(t)-P^{-1}V(t)P^{-1}\\
&=C^\top(t)Q(t)C(t)-P^{-1}V(t)P^{-1}-P^{-1}(A-\\
&\quad PC^\top(t)Q(t)C(t))-(A-PC^\top(t)Q(t)C(t))^\top P^{-1}.
\end{aligned}
\end{equation}
Now, consider the Lyapunov function candidate
\begin{align}
\label{eq:V_zeta}
\mathbf{V}(t,\zeta):=\frac{1}{\gamma}\zeta^\top P^{-1}(t)\zeta.
 \end{align}
It follows, from \eqref{eq:dzeta}, \eqref{ineq:g} and \eqref{eq:dP-1}, that the time-derivative of $\mathbf{V}(t,\zeta)$ in \eqref{eq:V_zeta} satisfies
\begin{equation}
\begin{aligned}
\dot{\mathbf{V}}(t,\zeta)
&=-\zeta^\top(C^\top(t)Q(t)C(t)+P^{-1}(t)V(t)P^{-1}(t))\zeta+\\
&\quad\frac{2}{\gamma^4}\zeta^\top P^{-1}B_2g(t,\tilde R,\tilde b_\omega)\\
&\leq-\frac{v_m}{\beta_2^2}\|\zeta\|^2+\frac{2c_g}{\beta_1\gamma^4}\|\zeta\|\\
&\leq-\frac{v_m}{2\beta_2^2}\|\zeta\|^2,\quad\forall\|\zeta\|\geq\frac{4c_g\beta_2^2}{\gamma^4\beta_1v_m}\\
&\leq-\frac{\gamma \beta_1v_m}{2\beta_2^2}\mathbf{V}(t,\zeta),\quad\forall\|\zeta\|\geq\frac{4c_g\beta_2^2}{\gamma^4\beta_1v_m}.
\label{ineq:dV_zeta}
\end{aligned}
\end{equation}
Let $c_\zeta$ and $T$ be any two positive constant scalars. Consider the set $\Omega=\{(t,\zeta): \mathbf{V}(t,\zeta)\leq \gamma^{-7}c_\zeta^2\beta_2^{-1}\}$. For any $(t,\zeta)\in\Omega$ one has
$
\|\zeta\|^2\leq\gamma\beta_2\mathbf{V}(t,\zeta)\leq ( c_\zeta\gamma^{-3})^2.
$
On the other hand, if $(t,\zeta)\notin\Omega$ and if we pick $\gamma\geq 4\beta_1^{-\frac{3}{2}}\beta_2^{\frac{5}{2}}c_gc_\zeta^{-1}v_m^{-1}$ then
\begin{align}\label{ineq:zeta1}
\|\zeta\|^2\geq\gamma\beta_1\mathbf{V}(\zeta)>\gamma^{-6}c_\zeta^2\beta_1\beta_2^{-1}\geq \left(\frac{4c_g\beta_2^2}{\gamma^4\beta_1v_m}\right)^2.
\end{align}
It follows in view of \eqref{ineq:dV_zeta} that for all $(t,\zeta)\notin\Omega$ one has $\dot{\mathbf{V}}(t,\zeta)\leq-\frac{\gamma \beta_1v_m}{2\beta_2^2}\mathbf{V}(t,\zeta)$. Hence, $(t,\zeta)$ must enter $\Omega$ before the following time:
\begin{align}\label{eq:T*}
T^*=\frac{2\beta_2^2}{\gamma \beta_1v_m}\ln\left(\frac{\gamma^7\beta_2\mathbf{V}(0,\zeta(0))}{c_\zeta^2}\right)
\end{align}
which can be tuned arbitrary small by increasing the value of $\gamma$. The following result immediately follows:
\begin{multline}\label{fact:bound:zeta}
\forall c_\zeta, T>0,\forall \zeta(0), \exists\gamma_1\geq 1\;\textrm{s.th.}\; \gamma\geq\gamma_1\Rightarrow\\ \|\zeta(t)\|\leq\gamma^{-3}c_\zeta,\;\forall t\geq T.
\end{multline}
Now we show that the gains can be tuned to guarantee forward invariance of the set $\mho(\epsilon)$. Let $\epsilon\in(\frac{1}{2},1)$ and let the initial conditions be such that $\tilde R(0)\in\mho(\epsilon)$. The time derivative of $|\tilde R|^2$, in view of \eqref{eq:dR} and \eqref{eq:dRhat} and making use of \cite[Lemma 2]{Berkane2017HybridSO3}, satisfies
\begin{equation}
\begin{aligned}
\nonumber
\frac{d}{dt}|\tilde R|^2&=-\frac{1}{4}\mathbf{tr}(\tilde R[-\hat R(\tilde b_\omega+k_1\sigma_2)]_\times)\\
&=-\frac{1}{4}\mathbf{tr}(\mathbf{P}_{\mathfrak{so}(3)}(\tilde R)[-\hat R(\tilde b_\omega+k_1\sigma_2)]_\times)\\
\nonumber
&=-\frac{1}{2}\psi(\tilde R)^\top\hat R(\tilde b_\omega+k_1\sigma_2)\\
&\leq\|\tilde b_\omega\|+k_1\|\sigma_2\|\\
&\leq c_b+k_1(\rho_1\|m_I\|^2+\rho_2c_2\hat c_2):=c_R.
\end{aligned}
\end{equation}
Let us define $t_R:=(\epsilon^2-|\tilde R(0)|^2)/c_R$. Hence, for all $0\leq t\leq t_R$, one has $\tilde R(t)\in\mho(\epsilon)$. Pick  $0<T\leq t_R$ and $c_\zeta:=\min(\bar c_\zeta/k_1,\hat c_2-c_2)$ for some arbitrary $\bar c_\zeta>0$. By \eqref{fact:bound:zeta} there exists $\gamma^*\geq1$ such that if one chooses $\gamma\geq\gamma^*$ then $\|\zeta(t)\|\leq \gamma^{-3}c_\zeta$ for all $t\geq T$. In this case, one has
\begin{equation}
\begin{aligned}
\|B_2^\top\hat x\|&=\|\hat a_I\|
=\|a_I-\hat a_I-a_I\|=\|B_2^\top L_\gamma\zeta-a_I\|\\
&=\|\gamma^3B_2^\top\zeta-a_I\|\leq \gamma^3\|\zeta\|+c_2\leq c_\zeta+c_2\leq \hat c_2.
\end{aligned}
\end{equation}
Consequently, for all $t\geq T$, one has $\mathbf{sat}_{\hat c_2}(\hat a_I(t))=\hat a_I(t)$.  It follows that the innovation term $\sigma_2$ in \eqref{eq:sigma2:1}-\eqref{eq:sigma2:2} is written as follows:
\begin{equation}
\begin{aligned}
\sigma_2
&=\rho_1(m_B\times\hat R^\top m_I)+\rho_2(a_B\times\hat R^\top \hat a_I)\\
&=\rho_1(m_B\times\hat R^\top m_I)+\rho_2(a_B\times\hat R^\top a_I)+\\
&\qquad\rho_2(a_B\times\hat R^\top (\hat a_I- a_I)\\
\label{eq:navigation:sigmaR2}
&=2\hat R^\top\psi(M\tilde R)-\rho_2(a_B\times\hat R^\top B_2^\top L_\gamma\zeta),
\end{aligned}    
\end{equation}
where we defined $M:=\rho_1 m_Im_I^\top+\rho_2a_Ia_I^\top$ and used  \cite[Proposition 3]{Berkane2017ConstructionStabilization} to derive the last equation. Note that $M$ is positive semidefinite and has rank equals $2$ (by Assumption \ref{assumption::obsv}). It follows that
\begin{equation}
\begin{aligned}
&\frac{d}{dt}|\tilde R|^2
=-k_1\psi(\tilde R)^\top\psi(M\tilde R)-\frac{1}{2}\psi(\tilde R)^\top(\hat R\tilde b_\omega-\\
&\qquad k_1\rho_2(\hat Ra_B)\times B_2^\top L_\gamma\zeta)\\
&\leq-4k_1\lambda_{\min}^{\E(M)}|\tilde R|^2(1-|\tilde R|^2)+|\tilde R|\|\tilde b_\omega\|+\\
&\qquad+\gamma^3 k_1\rho_2c_2|\tilde R|\|\zeta\|\\
&\leq -4k_1\lambda_{\min}^{\E(M)}|\tilde R(t)|^2(1-|\tilde R(t)|^2)+c_b+\rho_2c_2\bar c_\zeta
\end{aligned}    
\end{equation}
where inequalities from \cite[Lemma 2]{Berkane2017HybridSO3} have been used with $\E(M):=\frac{1}{2}(\tr(M)-M^\top)$ for any $M$. Note that the matrix $\E(M)$ is positive definite in view of Assumption \ref{assumption::obsv}. Now assume that $|\tilde R(t)|=\epsilon$ and $k_1>(c_b+\rho_2c_2\bar c_\zeta)/(4\lambda_{\min}^{\E(M)}\epsilon^2(1-\epsilon^2))$ then, for all $t\geq T$, one has
\begin{align*}
\frac{d}{dt}|\tilde R(t)|^2&\leq-4k_1\lambda_{\min}^{A}\epsilon^2(1-\epsilon^2)+c_b+\rho_2c_2\bar c_\zeta<0.
\end{align*}
This implies that $|\tilde R(t)|$ is strictly decreasing whenever $|\tilde R(t)|=\epsilon$. It follows from the continuity of the solutions that $\tilde R(t)$ will not leave the ball $\mho(\epsilon)$ for all $t\geq T$. Recall also that $|\tilde R(t)|\leq\epsilon$ for all $t\leq T$ (since $T\leq t_R$). This implies that the set $\mho(\epsilon)$ is forward invariant. Now, let us show exponential convergence. Consider the following Lyapunov function candidate:
\begin{multline}
\label{W}
\mathbf{W}(\zeta,\tilde R,\hat R,\tilde b_\omega):=|\tilde R|^2+\frac{\mu_1k_1}{2k_2}\tilde b_\omega^\top\tilde b_\omega+\mu_1\tilde b_\omega^\top\hat R^\top\psi(\tilde R)+\\\gamma^7\mathbf{V}(t,\zeta),
\end{multline}
where $\mu_1$ is some positive constant scalar and $\mathbf{V}(t,\zeta)$ is defined in \eqref{eq:V_zeta}.  Using the fact that $\|\psi(\tilde R)\|\leq 2|\tilde R|$, it can be verified that $\mathbf{W}$ satisfies the quadratic inequality $\varsigma^\top P_1\varsigma\leq \mathbf{W}\leq \varsigma^\top P_2 \varsigma$ where the matrices $P_1$ and $P_2$ are given by
\begin{align*}
P_1=\begin{bmatrix}
1&-\mu_1&0\\
-\mu_1&\frac{\mu_1 k_1}{2k_2}&0\\
0&0&\frac{\gamma^6}{\beta_2}
\end{bmatrix},
P_2=\begin{bmatrix}
1&\mu_1&0\\
\mu_1&\frac{\mu_1 k_1}{2k_2}&0\\
0&0&\frac{\gamma^6}{\beta_1}
\end{bmatrix}.
\end{align*}
Let us compute the time derivative of the cross term $\mathfrak{X}:=\tilde b_\omega^\top\hat R^\top\psi(\tilde R)$. First, one has
\begin{align*}\nonumber
				\dot{\mathfrak{X}}
                &=\tilde b_\omega^\top\hat R^\top \E(\tilde R)\left(-\hat R\tilde b_\omega-k_1\hat R\sigma_2)\right)-\\
                &\tilde b_\omega^\top[\omega+\tilde b_\omega+k_1\sigma_2]_\times\hat R^\top\psi(\tilde R)+\\
                &~\mathbf{P}_{c_5}^{\epsilon_b}\big(\hat b_\omega,k_2\sigma_2\big)^\top\hat R^\top\psi(\tilde R)\nonumber.
\end{align*}
In addition, using \cite[(31), Lemma 2]{Berkane2017HybridSO3}, one has the following bound
\begin{align*}
-\tilde b_\omega^\top\hat R^\top \E(\tilde R)\hat R\tilde b_\omega
&=-\|\tilde b_\omega\|^2+\tilde b_\omega^\top\hat R^\top (I-\E(\tilde R))\hat R\tilde b_\omega\\
&\leq-\|\tilde b_\omega\|^2+2c_b^2|\tilde R|^2.
\end{align*}
Moreover, in view of \eqref{eq:navigation:sigmaR2}, the following upper bound for $\sigma_2$ can be derived
\begin{align}
\|\sigma_2\|\leq 4\lambda_{\max}^{\E(M)}|\tilde R|+\rho_2 c_2\gamma^3\|\zeta\|.
\end{align}
It follows, using \cite[(18)]{Berkane2017AttitudeMeasurementsb}, that  
\begin{align*}
&-k_1\tilde b_\omega^\top\hat R^\top \E(\tilde R)\hat R\sigma_2=k_1\tilde b_\omega^\top\hat R^\top (I-\E(\tilde R))\hat R\sigma_2-\\
&k_1\tilde  b_\omega^\top\sigma_2\leq k_1|\tilde R|^2\tilde  b_\omega^\top\sigma_2+k_1\sqrt{2}|\tilde R|\|\tilde b_\omega\|\|\sigma_2\|-\\
&k_1\tilde  b_\omega^\top\sigma_2\leq-k_1\tilde  b_\omega^\top\sigma_2+8k_1\lambda_{\max}^{\E(M)}c_b(\sqrt{2}+2)|\tilde R|^2+\\
&\qquad 2k_1\gamma^3\rho_2c_bc_2(\sqrt{2}+2)|\tilde R|\|\zeta\|.
\end{align*}
Besides, the following bounds are easily derived
\begin{align*}
&-\tilde b_\omega^\top[\omega+\tilde b_\omega+k_1\sigma_2]_\times\hat R^\top\psi(\tilde R)\\
&=-\tilde b_\omega^\top[\omega+k_1\sigma_2]_\times\hat R^\top\psi(\tilde R)\\
&\leq c_\omega\|\tilde b_\omega\|\|\psi(\tilde R)\|+k_1\|\tilde b_\omega\|\|\psi(\tilde R)\|\|\sigma_2\|\\
&\leq 2c_\omega\|\tilde b_\omega\||\tilde R|+16k_1\lambda_{\max}^{\E(M)}c_b
|\tilde R|^2+\\
&\qquad+4k_1\gamma^3\rho_2c_bc_2|\tilde R|\|\zeta\|.
\end{align*}
and
\begin{align*}
&\mathbf{P}_{c_5}^{\epsilon_b}\big(\hat b_\omega,k_2\sigma_2\big)^\top\hat R^\top\psi(\tilde R)\leq k_2\|\sigma_2\|\|\psi(\tilde R)\|\\
&\leq 4k_2\lambda_{\max}^{\E(M)}|\tilde R|^2+2k_2\gamma^3\rho_2c_2|\tilde R|\|\zeta\|.
\end{align*}
Consequently, one deduces that
\begin{multline}
\label{ineq:observer2:Xi}
\dot{\mathfrak{X}}\leq-\|\tilde b_\omega\|^2-k_1\tilde  b_\omega^\top\sigma_2+(\alpha_1+k_1\alpha_2)|\tilde R|^2+\\
\gamma^3(\alpha_3+k_1\alpha_4)|\tilde R|\|\zeta\|+2c_\omega\|\tilde b_\omega\||\tilde R|,
\end{multline}
where $\alpha_1=8c_b^2+4k_2\lambda_{\max}^{\E(M)}$, $\alpha_2=8\lambda_{\max}^{\E(M)}c_b(\sqrt{2}+4))$, $\alpha_3=2k_2\rho_2c_2$ and $\alpha_4=2\rho_2c_bc_2(\sqrt{2}+4)$. Consequently, in view of the above obtained results, one has
\begin{equation}
\begin{aligned}
&\dot{\mathbf{W}}\leq-4k_1\lambda_{\min}^{\E(M)}(1-\epsilon^2)|\tilde R|^2-\mu_1\|\tilde b_\omega\|^2\\
&+2\gamma^3 c_2\beta_1^{-1}\|\zeta\|\|\tilde b_\omega\|
+(1+2\mu_1c_\omega)|\tilde R|\|\tilde b_\omega\|\\
&+\mu_1(\alpha_1+k_1\alpha_2)|\tilde R|^2-\gamma^7v_m\beta_2^{-2}\|\zeta\|^2\\
&+\gamma^3(k_1\rho_2c_2+4\sqrt{2}\beta_1^{-1}c_3+\mu_1(\alpha_3+k_1\alpha_4))|\tilde R|\|\zeta\|\\
&=-\varsigma_{12}^\top P_{12}\varsigma_{12}-\varsigma_{13}^\top P_{13}\varsigma_{13}-\varsigma_{23}^\top P_{23}\varsigma_{23}
\end{aligned}    
\end{equation}
where $\varsigma_{ij}=[\varsigma_i,\varsigma_j]^\top$ and the matrices $P_{ij}$ are given by
{\small
\begin{align*}
P_{12}&=\begin{bmatrix}
k_1\big(2\lambda_{\min}^{\E(M)}(1-\epsilon^2)-\mu_1\alpha_2\big)-\mu_1\alpha_1&\>-(\frac{1}{2}+c_\omega\mu_1)\\
*&\frac{\mu_1}{2}
\end{bmatrix}\\
P_{13}&=\begin{bmatrix}
2k_1\lambda_{\min}^{\E(M)}(1-\epsilon^2)&\quad -\frac{\gamma^3}{2}(k_1\rho_2c_2+4\sqrt{2}\beta_1^{-1}c_3+\\
&\qquad \qquad +\mu_1(\alpha_3+k_1\alpha_4))
\\
*&\frac{\gamma^7v_m}{2\beta_2^2}
\end{bmatrix}\\
P_{23}&=\begin{bmatrix}
\frac{\mu_1}{2}&-\gamma^3\beta_1^{-1}c_2\\
-\gamma^2\beta_1{-1}c_2&\frac{\gamma^7v_m}{2\beta_2^2}
\end{bmatrix}.
\end{align*}
}%
Now, if we pick $\mu_1>0$ such that $\mu_1<\lambda_{\min}^{\E(M)}(1-\epsilon^2)/\alpha_2$ and choose the gains $k_1$ and $\gamma$ such that
{\small
\begin{align*}
k_1&>\max\left\{2\mu_1 k_2,\frac{2\alpha_1\mu_1^2+(1+2c_\omega\mu_1)^2}{2\mu_1\lambda_{\min}^{\E(M)}(1-\epsilon^2)}\right\},\\
\gamma&>\max\Bigg\{ \frac{4\beta_2^2c_2^2}{\mu_1v_m\beta_1^2},\\
&\qquad \qquad \qquad \frac{\beta_2^2(k_1\rho_2c_2+4\sqrt{2}\beta_1^{-1}c_3+\mu_1(\alpha_3+k_1\alpha_4)^2}{4k_1\lambda_{\min}^{\E(M)}(1-\epsilon^2)v_m}\Bigg\},
\end{align*}}
then we can verify that matrices $P_1, P_2, P_{12}, P_{13}$ and $P_{23}$ are all positive definite. The exponential stability immediately follows.

\section{Proof of Lemma \ref{lemma:obsv:general}}\label{appendix:lemma:obsv:general}
Sufficiency: first note that 
$C^\top(t)C(t)=H^\top C_p^\top(t)C_p(t)H$ with $H=[I_3,0_{3\times 3},0_{3\times 3}]$. Since $A$ is constant, the observability Gramian of the pair $(A,C(\cdot))$ can be written as 
\begin{multline*}
W(t,t+\delta)=\\
=\int_t^{t+\delta}\exp(A(s-t))^\top H^\top C_p^\top(s)C_p(s)H\exp(A(s-t))ds.
\end{multline*}
Moreover, since $A$ is nilpotent with index $3$ ({\it i.e.,} $A^3=0$), the observability matrix of the pair $(A,H)$ satisfies $\textrm{rank}(H,HA,HA^2)=\textrm{rank}(I_9)=9$ and, hence, the pair $(A,H)$ is Kalman observable. Finally, all the eigenvalues of $A$ are zero (thus real). It follows, by direct application of \cite[Lemma 2.7]{Hamel2017PositionMeasurements}, that the PE condition \eqref{condition:pe:general} guarantees that the pair $(A,C(\cdot))$ is uniformly observable. \\
Necessity: we proceed by contradiction and assume that the PE condition of Lemma \ref{lemma:obsv:general} does not hold, {\it i.e.,} 
\begin{align}
    \forall\delta,\mu>0,\exists t\geq 0:\min_{z\in\mathbb{S}^2}\int_{t}^{t+\delta}\|C_p(s)z\|^2ds<\mu.
\end{align}
Consider a sequence $\{\mu_q\}_{q\in \mathbb{N}}$ of positive numbers converging to zero, and an arbitrary positive scalar $\delta$. Then, there must exist a sequence of time instants $\{t_q\}_{q\in \mathbb{N}}$ and a sequence of  unit vectors $\{z_q\}_{q\in \mathbb{N}}\subset\mathbb{S}^2$ such that  $\textstyle\frac{1}{\delta}  \int_{t_q}^{t_q+\delta} \|C(s)z_q\|^2 ds < \mu_q$  for any $q\in \mathbb{N}$. By the compactness of $\mathbb{S}^2$, there exists a sub-sequence of $\{z_q\}_{q\in \mathbb{N}}$ which converges to a limit $z\in \mathbb{S}^2$. Moreover, since   $C_p(t)$ is bounded and the interval of integration is fixed,  it follows that
\begin{align}
    \lim_{p\to\infty}\int_{t_p}^{t_p+\delta}\|C_p(s)z\|^2ds=0.
\end{align}
Let $\epsilon,\delta>0$ be arbitrary. Let $p$ be large enough such that 
$\textstyle\int_{t_p}^{t_p+\delta}\|C_p(s)z\|^2ds\leq\epsilon$. Now pick $x=[z^\top,0_{1\times 3},0_{1\times 3}]^\top\neq 0$. It follows that $A^kx=0$ for all $k\geq 1$ and $H\exp(At)x=Hx=z$.
Therefore
\begin{multline*}
    \int_{t_p}^{t_p+\delta}\|C_p(s)H\exp(A(s-t))x\|^2ds=\\\int_{t_p}^{t_p+\delta}\|C_p(s)z\|^2ds\leq\epsilon,
\end{multline*}
which shows that $(A,C(\cdot))$ is not uniformly observable.

\section{Proof of Lemma \ref{lemma:obsv}}\label{appendix:lemma:obsv}
Let us prove the result by contradiction.  Assume that \eqref{condition:pe:bearings} is not satisfied. That is to say,
\begin{align}\label{eq:not-pe1}
    \forall\delta,\mu>0,\exists t\geq 0,\frac{1}{\delta}\int_t^{t+\delta}\sum_{i=1}^n\Pi(R_iy_i(s))ds+\alpha e_3e_3^\top<\mu I_3.
\end{align}
Let $\{\mu_p\}_{p\in\mathbb{N}}$ be a sequence of positive numbers converging to zero and let $\delta>0$ be arbitrary. In view of \eqref{eq:not-pe1}, there must exist a sequence of times $\{t_p\}_{p\in\mathbb{N}}$ and a sequence of unit vectors $\{z_p\}_{p\in\mathbb{N}}$ such that 
\begin{align}
    \forall p\in\mathbb{N},\quad\frac{1}{\delta}\int_{t_p}^{t_p+\delta}\sum_{i=1}^n\|\Pi(R_iy_i(s))z_p\|^2ds+\alpha (e_3^\top z_p)^2<\mu_p.
\end{align}
Since $\mathbb{S}^2$ is compact, there must exit a sub-sequence of $\{z_p\}_{p\in\mathbb{N}}$ that converges to some limit unit vector $\bar z\in\mathbb{S}^2$. It follows that
\begin{align}
    \lim_{p\to\infty}\frac{1}{\delta}\int_{t_p}^{t_p+\delta}\sum_{i=1}^n\|\Pi(R_iy_i(s))\bar z\|^2ds+\alpha (e_3^\top\bar z)^2=0.
\end{align}
This is equivalent to 
\begin{align}
\label{eq:not-pe3}
&\alpha (e_3^\top\bar z)^2=0\\
\label{eq:not-pe2}
&\frac{1}{\delta}\lim_{p\to\infty}\int_0^{\delta}\|\Pi(R_iy_i(t_p+s))\bar z\|^2ds=0, i=1,\cdots,n.
\end{align}

Now, since $v$ is bounded, the function $\|\Pi(R_iy_i(t_p+s))\bar z\|^2$ is uniformly continuous and thus Cauchy-continuous. This implies that $\|\Pi(R_iy_i(t_p+s))\bar z\|^2$ is a Cauchy sequence of continuous functions and $ \lim_{p\to\infty}\|\Pi(R_iy_i(t_p+s))\bar z\|^2$ exists. Applying Lebesgue theorem, one has 
	$\textstyle
	\lim_{p\to\infty}\int_0^{\delta}\|\Pi(R_iy_i(t_p+s))\bar z\|^2ds   =  \int_0^{\delta}\lim_{p\to\infty}\|\Pi(R_iy_i(t_p+s))\bar z\|^2ds =0. $
	Since the function $\lim_{p\to\infty}\|\Pi(R_iy_i(t_p+s))\bar z\|^2$ is uniformly continuous and non-negative, it follows from
\eqref{eq:not-pe2} that
\begin{align}\label{eq:not-pe4}
    \lim_{p\to\infty}\|\Pi(R_iy_i(t_p+s))\bar z\|^2=0,\quad\forall s\in(0,\delta), i=1,\cdots
\end{align}
Now, since $\|\Pi(y)x\|^2=x^\top\Pi(y)x=-x^\top[y]^2_{\times}x=\|y\times x\|^2$ for all $x,y\in\mathbb{S}^2$, it follows that 
\begin{align}
    \lim_{p\to\infty}\|R_iy_i(t_p+s)\times\bar z\|=0,\quad\forall s\in(0,\delta), i=1,\cdots
\end{align}
This also implies that 
\begin{multline}
    \forall\zeta>0,\exists p^*,\forall p\geq p^*, \|R_iy_i(t_p+s)\times\bar z\|<\zeta,\\\forall s\in(0,\delta), i=1,\cdots
\end{multline}
Let $s\in(0,\delta/2)$ and pick $\zeta=\epsilon/2$. Then, there exists $p^*$ such that for all $p\geq p^*$
\begin{align}\label{eq:not-pe5}
    \|R_iy_i(t_p+s)\times\bar z\|<\epsilon/2,\\
    \label{eq:not-pe6}
    \|R_jy_j(t_p+\delta-s)\times\bar z\|<\epsilon/2,
\end{align}
for any $i,j$. Now since $\|x\times y\|^2=-y^\top[x]_{\times}^2y=y^\top (I_3-xx^\top)y=1-(x^\top y)^2$, for all $x,y\in\mathbb{S}^2$, and using the result in \citep{wang1994trace}, we obtain
\begin{equation}
    \|x\times y\|\leq\|x\times z\|+\|z\times y\|, \forall x,y,z\in\mathbb{S}^2
\end{equation}
Using this latter fact, inequalities \eqref{eq:not-pe5}-\eqref{eq:not-pe6} imply that
\begin{align}
    \|R_iy_i(t_p+s)\times R_jy_j(t_p+\delta-s)\|<\epsilon.
\end{align}
Since $\delta$ can be arbitrary large and $s$ can be arbitrary small, this last equation contradicts item i) of the Lemma. Furthermore, note that for any $t\geq 0$ and any $i,j$ we can write
\begin{align*}
    &\|\bar z\times (p_i-p_j)\|^2
    = \|\bar z\times (p(t)-p_i)\|^2+ \|\bar z\times (p(t)-p_j)\|^2\\
    &\qquad-2(\bar z\times (p(t)-p_i))^\top(\bar z\times (p(t)-p_j))\\
    &=\beta_i^2(t)\|\bar z\times R_iy_i(t) \|^2+\beta_j^2(t)\|\bar z\times R_jy_j(t)\|^2\\
    &-2\beta_i(t)\beta_j(t)(\bar z\times R_iy_i(t))^\top(\bar z\times R_jy_j(t))
\end{align*}
where $\beta_i(t)=\|p(t)-p_i\|$ and $\beta_j(t)=\|p(t)-p_j\|$. However, $\lim_{p\to\infty}\bar z\times R_iy_i(t_p+s)=0$ for $i=1,\cdots$. By selecting $t=t_p+s$ and letting $p$ go to infinity in the above equation, it follows that $\|\bar z\times (p_i-p_j)\|^2=0$ and, thus, $\bar z$ is parallel to $(p_i-p_j)$ for any $i$ and $j$. If we have three source points $p_i, p_j$ and $p_k$ that are not aligned, it is not possible to have $\bar z$ to be parallel to $(p_i-p_j)$ and $(p_i-p_k)$ simultaneously. Therefore, item ii) of the lemma holds. \\
If the altimeter is used, \eqref{eq:not-pe3} implies that $e_3^\top\bar z=e_3^\top (p_i-p_j)=0$ which contradicts item iv) of the lemma. Finally,  in view of \eqref{eq:not-pe4}, we have $\lim_{p\to\infty}((R_iy_i(t_p+s))^\top\bar z)^2=1$ and $\lim_{p\to\infty} \Pi(R_iy_i(t_p+s))\bar z=0$. Therefore,
\begin{align}
    &\lim_{p\to\infty} e_3^\top \Pi(R_iy_i(t_p+s))\bar z\\
    &=e_3^\top\bar z-\lim_{p\to\infty}(e_3^\top R_iy_i(t_p+s))(R_iy_i(t_p+s)^\top\bar z)
    \\& =\pm \lim_{p\to\infty}(e_3^\top R_iy_i(t_p+s))\\
    &=0
\end{align}
for all $s\in(0,\delta)$, where we have used the fact that $e_3^\top\bar z=0$ if the altimeter is used.  This contradicts item iii) of the lemma. 
\bibliographystyle{apalike}
\bibliography{references2.bib}

\begin{thebibliography}{}

\bibitem[Batista et~al., 2017]{Batista2017RelaxedRealizations}
Batista, P., Petit, N., Silvestre, C., and Oliveira, P. (2017).
\newblock {Relaxed conditions for uniform complete observability and
  controllability of LTV systems with bounded realizations}.
\newblock {\em IFAC-PapersOnLine}, 50(1):3598--3605.

\bibitem[Batista et~al., 2011]{Batista2011SingleDesign}
Batista, P., Silvestre, C., and Oliveira, P. (2011).
\newblock {Single range aided navigation and source localization: Observability
  and filter design}.
\newblock {\em Systems and Control Letters}, 60(8):665--673.

\bibitem[Batista et~al., 2012]{Batista2012GloballyEstimation}
Batista, P., Silvestre, C., and Oliveira, P. (2012).
\newblock {Globally exponentially stable cascade observers for attitude
  estimation}.
\newblock {\em Control Engineering Practice}, 20(2):148--155.

\bibitem[Batista et~al., 2013]{Batista2013GloballyMeasurements}
Batista, P., Silvestre, C., and Oliveira, P. (2013).
\newblock {Globally exponentially stable filters for source localization and
  navigation aided by direction measurements}.
\newblock {\em Systems {\&} Control Letters}, 62(11):1065--1072.

\bibitem[Batista et~al., 2015]{Batista2015NavigationMeasurements}
Batista, P., Silvestre, C., and Oliveira, P. (2015).
\newblock {Navigation systems based on multiple bearing measurements}.
\newblock {\em IEEE Transactions on Aerospace and Electronic Systems},
  51(4):2887--2899.

\bibitem[Batista et~al., 2016]{Batista2016TightlySystem}
Batista, P., Silvestre, C., and Oliveira, P. (2016).
\newblock {Tightly coupled long baseline/ultra-short baseline integrated
  navigation system}.
\newblock {\em International Journal of Systems Science}, 47(8):1837–1855.

\bibitem[Berkane and Abdelhamid, 2019]{Berkane2017AttitudeMeasurementsb}
Berkane, S. and Abdelhamid, T. (2019).
\newblock {Attitude Estimation with Intermittent Measurements}.
\newblock {\em Automatica}, 105:415--421.

\bibitem[Berkane et~al., 2016]{Berkane2016GlobalSystems}
Berkane, S., Abdessameud, A., and Tayebi, A. (2016).
\newblock {Global exponential angular velocity observer for rigid body
  systems}.
\newblock In {\em IEEE 55th Conference on Decision and Control}, pages
  4154--4159.

\bibitem[Berkane et~al., 2017]{Berkane2017HybridSO3}
Berkane, S., Abdessameud, A., and Tayebi, A. (2017).
\newblock {Hybrid Attitude and Gyro-bias Observer Design on SO(3)}.
\newblock {\em IEEE Transactions on Automatic Control}, 62(11):6044--6050.

\bibitem[Berkane and Tayebi, 2017a]{Berkane2017AttitudeMeasurements}
Berkane, S. and Tayebi, A. (2017a).
\newblock {Attitude and Gyro Bias Estimation Using GPS and IMU Measurements}.
\newblock In {\em IEEE 56th Conference on Decision and Control}, pages
  2402--2407.

\bibitem[Berkane and Tayebi, 2017b]{Berkane2017ConstructionStabilization}
Berkane, S. and Tayebi, A. (2017b).
\newblock {Construction of Synergistic Potential Functions on SO(3) with
  Application to Velocity-Free Hybrid Attitude Stabilization}.
\newblock {\em IEEE Transactions on Automatic Control}, 62(1):495--501.

\bibitem[Berkane and Tayebi, 2017c]{Berkane2017OnSO3}
Berkane, S. and Tayebi, A. (2017c).
\newblock {On the Design of Attitude Complementary Filters on SO(3)}.
\newblock {\em IEEE Transactions on Automatic Control}, 63(3):880--887.

\bibitem[Berkane and Tayebi, 2019]{berkane2019position}
Berkane, S. and Tayebi, A. (2019).
\newblock Position, velocity, attitude and gyro-bias estimation from imu and
  position information.
\newblock In {\em 18th European Control Conference (ECC)}, pages 4028--4033.

\bibitem[Bonnabel et~al., 2006]{Bonnabel2006}
Bonnabel, S., Martin, P., and Rouchon, P. (2006).
\newblock {A non-linear symmetry-preserving observer for velocity-aided
  inertial navigation}.
\newblock In {\em Proceedings of the American Control Conference}, pages
  2910--2914.

\bibitem[Bristeau et~al., 2010]{Bristeau2010DesignObservability}
Bristeau, P.~J., Petit, N., and Praly, L. (2010).
\newblock {Design of a navigation filter by analysis of local observability}.
\newblock {\em Proceedings of the IEEE Conference on Decision and Control},
  pages 1298--1305.

\bibitem[Bryne et~al., 2017]{Bryne2017NonlinearAspects}
Bryne, T.~H., Hansen, J.~M., Rogne, R.~H., Sokolova, N., Fossen, T.~I., and
  Johansen, T.~A. (2017).
\newblock {Nonlinear Observers for Integrated Implementation Aspects}.
\newblock {\em IEEE Control Systems Magazine}, 37(3):59--86.

\bibitem[Bucy, 1972]{Bucy1972TheBounds}
Bucy, R.~S. (1972).
\newblock {The riccati equation and its bounds}.
\newblock {\em Journal of Computer and System Sciences}, 6(4):343--353.

\bibitem[Crassidis, 2006]{Crassidis2006Sigma-pointNavigation}
Crassidis, J.~L. (2006).
\newblock {Sigma-point Kalman filtering for integrated GPS and inertial
  navigation}.
\newblock {\em IEEE Transactions on Aerospace and Electronic Systems},
  42(2):750--756.

\bibitem[Esfandiari and Khalil, 1992]{esfandiari1992output}
Esfandiari, F. and Khalil, H.~K. (1992).
\newblock Output feedback stabilization of fully linearizable systems.
\newblock {\em International Journal of control}, 56(5):1007--1037.

\bibitem[Farrell, 2008]{Farrell2008AidedSensors}
Farrell, J. (2008).
\newblock {\em {Aided navigation: GPS with high rate sensors}}.
\newblock McGraw-Hill, Inc.

\bibitem[Grip et~al., 2013]{Grip2013}
Grip, H.~F., Fossen, T.~I., Johansen, T.~A., and Saberi, A. (2013).
\newblock {Nonlinear Observer for GNSS-Aided Inertial Navigation with
  Quaternion-Based Attitude Estimation}.
\newblock In {\em American Control Conference}, pages 272--279.

\bibitem[Gryte et~al., 2017]{Gryte2017RobustGPS}
Gryte, K., Hansen, J.~M., Johansen, T., and Fossen, T.~I. (2017).
\newblock {Robust navigation of UAV using inertial sensors aided by UWB and RTK
  GPS}.
\newblock In {\em AIAA Guidance, Navigation, and Control Conference}, pages
  1--16.

\bibitem[Hamel and Samson, 2017]{Hamel2017PositionMeasurements}
Hamel, T. and Samson, C. (2017).
\newblock {Position estimation from direction or range measurements}.
\newblock {\em Automatica}, 82:137--144.

\bibitem[Hamel and Samson, 2018]{Hamel2018RiccatiProblem}
Hamel, T. and Samson, C. (2018).
\newblock {Riccati Observers for the Nonstationary PnP Problem}.
\newblock {\em IEEE Transactions on Automatic Control}, 63(3):726--741.

\bibitem[Hamer and D'Andrea, 2018]{Hamer2018Self-calibratingLocalization}
Hamer, M. and D'Andrea, R. (2018).
\newblock {Self-calibrating ultra-wideband network supporting multi-robot
  localization}.
\newblock {\em IEEE Access}, 6:22292--22304.

\bibitem[Hansen et~al., 2018]{hansen2018nonlinear}
Hansen, J.~M., Johansen, T.~A., Sokolova, N., and Fossen, T.~I. (2018).
\newblock Nonlinear observer for tightly coupled integrated inertial navigation
  aided by rtk-gnss measurements.
\newblock {\em IEEE Transactions on Control Systems Technology},
  27(3):1084--1099.

\bibitem[Hua, 2010]{hua2010attitude}
Hua, M.-D. (2010).
\newblock Attitude estimation for accelerated vehicles using {GPS/INS}
  measurements.
\newblock {\em Control Engineering Practice}, 18(7):723--732.

\bibitem[Johansen and Fossen, 2015]{johansen2015nonlinear}
Johansen, T.~A. and Fossen, T.~I. (2015).
\newblock Nonlinear observer for inertial navigation aided by pseudo-range and
  range-rate measurements.
\newblock In {\em European Control Conference (ECC)}, pages 1673--1680.

\bibitem[Johansen et~al., 2017]{Johansen2017NonlinearMeasurements}
Johansen, T.~A., Hansen, J.~M., and Fossen, T.~I. (2017).
\newblock {Nonlinear Observer for Tightly Integrated Inertial Navigation Aided
  by Pseudo-Range Measurements}.
\newblock {\em Journal of Dynamic Systems, Measurement, and Control},
  139(1):1--10.

\bibitem[Kai et~al., 2017]{KAI2017}
Kai, J.-M., Allibert, G., Hua, M.-D., and Hamel, T. (2017).
\newblock Nonlinear feedback control of quadrotors exploiting first-order drag
  effects.
\newblock {\em IFAC-PapersOnLine}, 50(1):8189 -- 8195.
\newblock 20th IFAC World Congress.

\bibitem[Krstic et~al., 1995]{Krstic1995NonlinearDesign}
Krstic, M., Kanellakopoulos, I., and Kokotovic, P.~V. (1995).
\newblock {\em {Nonlinear and adaptive control design}}.
\newblock Wiley New York.

\bibitem[Mahony et~al., 2008]{Mahony2008NonlinearGroup}
Mahony, R., Hamel, T., and Pflimlin, J.-m. (2008).
\newblock {Nonlinear Complementary Filters on the Special Orthogonal Group}.
\newblock {\em IEEE Transactions on Automatic Control}, 53(5):1203--1218.

\bibitem[Morin et~al., 2017]{morin2017uniform}
Morin, P., Eudes, A., and Scandaroli, G. (2017).
\newblock Uniform observability of linear time-varying systems and application
  to robotics problems.
\newblock In {\em International Conference on Geometric Science of
  Information}, pages 336--344. Springer.

\bibitem[Roberts and Tayebi, 2011]{roberts2011}
Roberts, A. and Tayebi, A. (2011).
\newblock On the attitude estimation of accelerating rigid-bodies using {GPS}
  and {IMU} measurements.
\newblock In {\em 50th IEEE conference on decision and control and European
  Control conference}, pages 8088--8093.

\bibitem[Sabatini, 2006]{sabatini2006quaternion}
Sabatini, A.~M. (2006).
\newblock Quaternion-based extended kalman filter for determining orientation
  by inertial and magnetic sensing.
\newblock {\em IEEE transactions on Biomedical Engineering}, 53(7):1346--1356.

\bibitem[{Saberi} and {Sannuti}, 1990]{saberi1990}
{Saberi}, A. and {Sannuti}, P. (1990).
\newblock Observer design for loop transfer recovery and for uncertain
  dynamical systems.
\newblock {\em IEEE Transactions on Automatic Control}, 35(8):878--897.

\bibitem[Titterton et~al., 2004]{Titterton2004StrapdownTechnology}
Titterton, D., Weston, J.~L., and Weston, J. (2004).
\newblock {\em {Strapdown inertial navigation technology}}, volume~17.
\newblock IET.

\bibitem[Vik and Fossen, 2001]{vik2001nonlinear}
Vik, B. and Fossen, T.~I. (2001).
\newblock A nonlinear observer for gps and ins integration.
\newblock In {\em Proceedings of the 40th IEEE Conference on Decision and
  Control (Cat. No. 01CH37228)}, volume~3, pages 2956--2961.

\bibitem[Wang and Zhang, 1994]{wang1994trace}
Wang, B.-Y. and Zhang, F. (1994).
\newblock A trace inequality for unitary matrices.
\newblock {\em The American Mathematical Monthly}, 101(5):453--455.

\bibitem[Whitmore, 2000]{whitmore2000closed}
Whitmore, S.~A. (2000).
\newblock Closed-form integrator for the quaternion (euler angle) kinematics
  equations.
\newblock US Patent 6,061,611.

\bibitem[Whittaker and Crassidis, 2017]{Whittaker2017InertialRepresentations}
Whittaker, M. and Crassidis, J.~L. (2017).
\newblock {Inertial Navigation Employing Common Frame Error Representations}.
\newblock {\em AIAA Guidance, Navigation, and Control Conference}, pages 1--24.

\bibitem[Woodman, 2007]{Woodman2007AnNavigation}
Woodman, O.~J. (2007).
\newblock {An introduction to inertial navigation}.
\newblock {\em Technical Report, University of Cambridge Computer, Laboratory},
  (UCAM-CL-TR-696):1--37.

\bibitem[Zlotnik and Forbes, 2016]{Zlotnik2016NonlinearDirectly}
Zlotnik, D.~E. and Forbes, J.~R. (2016).
\newblock {Nonlinear Estimator Design on the Special Orthogonal Group using
  Vector Measurements Directly}.
\newblock {\em IEEE Transactions on Automatic Control}, 9286(c):1--12.

\end{thebibliography}
\end{document}